\documentclass{article}
\usepackage{amssymb,amsmath}
\usepackage[latin1]{inputenc}
\usepackage[T1]{fontenc}
\usepackage[french,english]{babel}
\usepackage{amsthm}
\usepackage{amsfonts}
\usepackage{pxfonts}
\author{Dominique MALICET}
\date{}
\title{Lyapunov exponent of random dynamical systems on the circle}
\newtheorem{thm1}{Theorem}
\newtheorem{lem1}{Lemma}[section]
\newtheorem{cor1}{Corollary}
\newtheorem{Def1}{Definition}[section]
\newtheorem{prop1}{Proposition}[section]
\newtheorem{rem}{Remark}[section]
\newcommand{\Z}{\mathbb{Z}}
\newcommand{\R}{\mathbb{R}}

\newcommand{\N}{\mathbb{N}}
\newcommand{\C}{\mathbb{C}}
\newcommand{\E}{\mathbb{E}}
\newcommand{\T}{\mathbb{T}}

\newcommand{\ep}{\varepsilon}

\newcommand{\pa}[1]{\left( #1 \right)}
\newcommand{\cro}[1]{\left[ #1 \right]}

\newenvironment{disarray}{\everymath{\displaystyle\everymath{}}\array}{\endarray}

\begin{document}
\maketitle
\begin{abstract}
We consider products of a i.i.d. sequence in a set $\{f_1,\ldots,f_m\}$ of preserving orientation diffeomorphisms of the circle. we can naturally associate a Lyapunov exponent $\lambda$. Under few assumptions, it is known that $\lambda\leq 0$ and that the equality holds if and only if $f_1,\ldots,f_m$ are simultaneously conjugated to rotations. In this paper, we state a quantitative version of this fact in the case where $f_1,\ldots,f_m$ are $C^k$ perturbations of rotations with rotation numbers $\rho(f_1),\ldots,\rho(f_m)$ satisfying a simultaneous diophantine condition in the sense of Moser \cite{Moser}: we give a precise estimate on $\lambda$ (Taylor expansion) and we prove that there exists a diffeomorphism $g$ and rotations $r_i$ such that $\mbox{dist}(gf_ig^{-1},r_i)\ll |\lambda|^{\frac{1}{2}}$ for $i=1,\ldots m$. We also state analog results for random products of matrices $2\times 2$, without diophantine condition.
\end{abstract}

\section{Statement of results}

\subsection{Lyapunov exponent of random product of diffeomorphisms of the torus}

 We consider the random compositions $g_n=f_{n-1}\circ\cdots \circ f_0$ where $(f_k)_{k\in\N}$ is a sequence of i.i.d. copies of some random diffeomorphism $f$ of the unidmensional torus $\T=\R/\Z$. The general expected behaviour under few assumptions is that alsmost surely, the random orbits $(g_n(x))_{n\in\N}$ distribute themselves toward a unique \textit{stationary probability measure} $\mu$ on $\T$, and that the derivatives $g_n'(x)$ decrease toward $0$ with a fixed exponential rate given by a \textit{Lyapunov exponent} $\lambda$ (we will recall the precise definitions). The objective is to estimate the measure $\mu$ and the number $\lambda$ when $f$ is the perturbation a random rotation, and to obtain by an explicit estimate that $\lambda$ is an obstruction to the existence of a linearization of $f$, that is to say a deterministic diffeomorphism $g$ such that $gfg^{-1}$ is a rotation.\\

Let us begin by introducing some notations: the circle is identified with the torus $\T=\R/\Z$. For $k\in\N$ we identify $C^k(\T)$ with the space of $1$-periodic $C^k$ maps from $\R$ into $\R$ endowed with its standard norm $\|\cdot\|_k$ defined by $\|\varphi\|_k=\sup_{j\leq k,x\in\R}|\phi^{(j)}(x)|$. In the same way $\mbox{Diff}_+^k(\T)$ is the space of increasing diffeomorphisms $f$ from $\R$ onto $\R$ on the form $f=Id+\varphi$ with $\varphi\in C^k(\T)$. Noting that the difference of two elements of $\mbox{Diff}_+^k(\T)$ belongs to $C^k(\T)$ allows to naturally endow $\mbox{Diff}_+^k(\T)$ with the metric $d_k$ defined by $d_k(f,g)=\|f-g\|_k$. With these definitions, a rotation of $\T$ of angle $\alpha$ is simply the translation $Id+\alpha$, that we denote $r_\alpha$.\\ 
 
A random diffeomorphism of $\T$ is a random variable valued in $\mbox{Diff}_+(\T)$. In the paper all the random variables are implicitely assumed defined on a same probability space $(\Omega,\mathcal{F},\mathbb{P})$. Let us recall the notions of stationary measure and Lyapunov exponent for a random diffeomorphism:

\begin{Def1}
Let $f$ be a random diffeomorphism of $\T$  valued in $\mbox{Diff}_+^k(\T)$ such that $\ln_+ \|f'\|_0\in L^1(\Omega)$. A probability measure $\mu$ on $\T$ is stationary for $f$ if $\E[f_*\mu]=\mu$ (such a measure always exists by Kakutani fixed point theorem). The associated (mean) Lyapunov exponent is
$$\lambda(\mu)=\E\int_\T\ln{|f'(x)|}d\mu(x).$$

\end{Def1}

We recall some known facts about stationary measures and Lyapunov exponents. We will not use them in this paper but it may enlighten the reader on their meaning and their interest.
\begin{prop1}\label{recap}
Let $f$ be a random diffeomorphism valued in $\mbox{Diff}_+^1(\T)$ such that $\ln_+ \|f'\|_0\in L^1(\Omega)$, and let $g_n=f_{n-1}\circ\cdots \circ f_0$, where $(f_k)_{k\in\N}$ is a sequence of i.i.d. copies of $f$.
\begin{itemize}
\item If $f$ is minimal in the sense that the unique closed sets of $\T$ almost surely invariant by $f$ are $\emptyset$ and $\T$, then the stationary measure is unique. (see \cite{Deroin-inter}, \cite{Malicet})
\item If there is a unique stationary measure $\mu$ for $f$ and so a unique Lyapunov exponent $\lambda=\lambda(\mu)$, then for every $x$ in $\T$ we have 
$$\frac{1}{n}\ln (g_n'(x))\xrightarrow[n\to +\infty]{}\lambda \mbox{ a.s.}$$
\item $\lambda(\mu)$ is a negative number unless maybe if almost every realization of $f$ preserves $\mu$ (it is an early version due to Crauel \cite{Crauel} of the so called ``invariance principle'' of Avila-Viana \cite{Avila}, both inspired by the linear version in the seminal paper \cite{Ledrappier} of Ledrappier).\\
 If $f$ is minimal, it implies the existence of a homeomorphism $h$ of $\T$ such that $hfh^{-1}$ is almost surely a rotation, and so implies in particular that a.e. realizations of $f$ commute.  
\end{itemize}
\end{prop1}
 
We are going to give an estimate for $\lambda(\mu)$ when $f$ is a perturbation of a random rotation. We need an arithmetical condition on the angle of the random rotation. We recall that a number $\alpha$ is diophantine if for some $A,\sigma>0$ we have $\mbox{dist}(q\alpha,\Z)\geq \frac{A}{|q|^{\sigma}}$ for any $q$ in $\Z-\{0\}$, definition generalized by Moser in \cite{Moser} where $m$ numbers $\alpha_1,\ldots,\alpha_m$ are said simultaneously diophantine if for some $A,\sigma>0$ we have $\sup_i \mbox{dist}(q\alpha_i,\Z)\geq \frac{A}{|q|^{\sigma}}$ for any $q$ in $\Z-\{0\}$ (in particular, it holds if at least one of the $\alpha_i$ is diophantine). Here we introduce a definition generalizing the classical notion of diophantine number for random variables.
\begin{Def1}\label{diophantien}
Let $\alpha$ be a random variable in $\T$. For any $A>0$ and $\sigma\geq 0$, we say that $\alpha$ is diophantine of type $(A,\sigma)$ if for any $q$ in $\Z-\{0\}$,
\begin{equation}\label{dio}\left\|\mbox{dist}(q\alpha,\Z)\right\|_{L^2(\Omega)}\geq \frac{A}{|q|^{\sigma}}.\end{equation}
We say that $\alpha$ is diophantine if there exists $A>0$ and $\sigma\geq 0$ such that $\alpha$ is diophantine of type $(A,\sigma)$.
\end{Def1}
\begin{rem}~

	\begin{itemize}
		\item If $\alpha$ is deterministic (i.e. is a constant random variable), then we obtain the classical definition of diophantine number, and if the set of realizations of $\alpha$ is a finite set $\{\alpha_1,\ldots,\alpha_m\}$, then $\alpha$ is diophantine if and only if $\alpha_1,\ldots,\alpha_m$ are simultaneously diophantine.
		\item If $\alpha$ has positive probability to be a diophantine number, then $\alpha$ is a diophantine random variable.
		\item At the contrary to the deterministic case, it can happen that $\sigma=0$. It is for exemple the case if $\alpha$ is uniform on $\T$ by a simple computation(or more generally if the law of $\alpha$ is not Lebesgue singular, by a consequence of Riemman Lebesgue Lemma)
	\end{itemize}
\end{rem}
To check the second point, consider the sets $E_{A,\sigma}$ of $x$ in $\T$ such that for every $q$ in $\Z^*, \mbox{dist}(qx,\Z)\geq\frac{A}{|q|^\sigma}$. If $\alpha$ has positive probability to be diophantine, then there must exist $A$ and $\sigma$ such that $\alpha$ belongs to $E_{A,\sigma}$ with positive probability $p$, and then: $\forall q\in \Z^*, \left\|\mbox{dist}(q\alpha,\Z)\right\|_{L^2(\Omega)}\geq \frac{A}{|q|^{\sigma}}\sqrt{p}$.\\

Our first theorem gives a precise estimate for the Lyapunov exponent of a random diffeormorphism $f=r_\alpha+\zetaup$ when $f$ is a perturbation (in a smooth sense) of order $\ep$ of a random rotation $r_\alpha$ with $\alpha$ diophantine. We obtain a quadratic estimate $\lambda=O(\ep^2)$ (instead of the obvious bound $\lambda=O(\ep)$) and a formula for the quadratic term. In the statement of the theorem, a term $O(M)$ means a term bounded by $CM$ with $C$ a constant depending only on $A$ and $\sigma$.
\begin{thm1}\label{lyapu}
Let $\alpha$ be a diophantine random variable of type $(A,\sigma)$. Then there exists an integer $k$ depending only on $\sigma$ such that for any random random diffeomorphism in $\mbox{Diff}_+^k(\T)$ on the form $f=r_\alpha+\zetaup$ and for any Lyapunov  exponent  $\lambda$ associated to any stationary measure of $f$, we have
	$$\lambda=-\frac{1}{2}\E\int_{\T}\left( \zetaup '+\etaup'-\etaup'\circ r_\alpha \right)^2dx+O(\ep^3)$$
	(and so $\lambda=O(\ep^2)$), where $\ep=\|d_k(f,r_\alpha)\|_{L^3(\Omega)}=\E\cro{d_k(f,r_\alpha)^3}^{\frac{1}{3}}$, and where $\etaup$ is a deterministic map depending linearly on $\zetaup$ and satisfying $|\eta '|=O(\ep)$. The non zero Fourier coefficients of $\etaup$ are given by the formula
	
	\begin{equation}\label{estilambda}
	\hat{\etaup}(p)=\frac{\E[\hat{\zetaup}(p)e^{-2i\pi p \alpha}]}{1-\E[e^{-2i\pi p\alpha}]}.	\end{equation}

\end{thm1}
The formula \ref{estilambda} can also be rewritten by Parseval identity as
$$\lambda=-\frac{1}{2}\E\sum_{p\in\Z^*}p^2\left|\hat{\zetaup}(p)+\frac{\E[\hat{\zetaup}(p)e^{-2i\pi p \alpha}]}{1-\E[e^{-2i\pi p\alpha}]}(1-e^{2i\pi\alpha})\right|^2+O(\ep^3)$$
\begin{rem}
	Our method can actually allow to obtain the higher terms in the Taylor expansion of $\lambda$, on the form $\lambda=\sum_{j=2}^{n-1} q_j(\zetaup)+O(\ep^n)$ where $q_j(\zetaup)$ is a $j$-linear form evaluated at $(\zetaup,\ldots,\zetaup)$.
\end{rem}

In the next theorem we prove that if $f$ is a random diffeomorphism close to rotations whose rotation number $\rho(f)$ is diophantine, then $\lambda$ measures in an explicit sense how much close to rotations $f$ can be (smoothly) conjugated by a deterministic diffeomorphism. Note that $\lambda$ is indeed a natural obstruction to the existence of such a diffeomorphism because $\lambda$ is invariant by conjugation.
\begin{thm1}\label{principal}
Let $(A,\sigma)$ be a couple of positive real numbers. There exists an integer $r$ depending only on $\sigma$ such that for any integer $K$ larger than $r$, there exists in $\mbox{Diff}_+^K(\T)$ a neighborhood $\mathcal{U}$ of the set of rotations such that for any random diffeomorphism $f$ valued in $\mathcal{U}$ whose rotation number $\alpha=\rho(f)$ is $(A,\sigma)$ diophantine, there exists in $\mbox{Diff}_+^{K-r}(\T)$ a (non random) diffeomorphism $h$ such that 
	$$\|d_0(hfh^{-1},r_\alpha)\|_{L^2(\Omega)}\leq 3|\lambda|^{\frac{1}{2}},$$
for any Lyapunov exponent $\lambda$ associated to a stationary measure of $f$, with $h$ satisfying $d_{K-r}(h,Id)\leq C \|d_K(f,r_\alpha)\|_{L^2(\Omega)}$ for some $C$ depending on $A$, $\sigma$ and $K$.
\end{thm1}
The constant $3$ in the inequality above is not optimal. By analyzing carefully our proof we could actually replace it by any number larger than $\sqrt{2}$. However the bound $|\lambda|^{\frac{1}{2}}$ is essentially optimal since by Theorem \ref{lyapu}, $|\lambda|^{\frac{1}{2}}=O(d_k(hfh^{-1},r_\alpha))$ for some integer $k$. The number $r$ represents the ``loss of derivative''. It can be explicited from our proof as an affin function of $\sigma$, though we did not try at all to obtain an optimal expression.
\begin{rem}
If $\lambda=0$ and $f$ is valued in a finite set $\{f_1,\ldots,f_m \}$ the theorem gives  a smooth diffeomorphism $h$ conjugating silmutaneously $f_1,\ldots,f_m$ to rotations. This particular case can actually be obtained by using a succession of already known results: $f$ is minimal by Denjoy theorem (the diophantine condition implies that at least  one of the rotation numbers $\rho(f_i)$ is irrational), so if $\lambda=0$  the maps $f_i$ are simultaneously $C^0$-conjugated to rotations $r_1,...,r_m$ and so pairwise commute (see Proposition \ref{recap}).Then one can use a result of Moser \cite{Moser} which generalizes the classical works of Arnold \cite{Arnold} and Moser on the linearization of a single map close to rotations in the case of several commuting maps, and which states that under the diophantine condition given in assumption, the conjugacy $h$ can be taken smooth and close to Identity with the estimate $d_{K-r}(h,Id)=O(\sup_j d_K(f_j,r_j))$.\end{rem}
Since the maps close to rotations almost commute, we can deduce from Theorem \ref{principal} the following corollary:
\begin{cor1}\label{coro}
Let $(A,\sigma)$ be a couple of positive real numbers.  Then there exists an integer $k$
 and a neighborhood $\mathcal{U}$ of the set of rotations in $\mbox{Diff}_+^{k}(\T)$ such that for any random diffeomorphism $f$ valued in $\mathcal{U}$, if $\alpha=\rho(f)$ is $(A,\sigma)$ diophantine then, by denoting by $\tilde{f}$ an independent copy of $f$ we have
	$$||d_{0}(f\circ \tilde{f},\tilde{f}\circ f)||_{L^2(\Omega)}\leq C|\lambda|^{\frac{1}{2}}$$
	for any Lyapunov exponent $\lambda$ associated to a stationary measure of $f$, where $C$ is a universal constant.
\end{cor1}
	
 By Theorem \ref{principal} there exists an integer $k$ and a neighborhood $\mathcal{U}$ of rotations in $\mbox{Diff}^{k}(\T)$ such that for $f$ valued in $\mathcal{U}$, there exists $h$ in $\mbox{Diff}_+^1(\T)$ with $\max(h',(h^{-1})')\leq 2$ such that $f_1=hfh^{-1}$ satisfies $\|d_0(f_1,r_\alpha)\|_{L^2(\Omega)}\leq 3|\lambda|^{\frac{1}{2}}$. Then, setting $\tilde{f_1}=h\tilde{f}h^{-1}$ and $\tilde{\alpha}=\rho(\tilde{f})$ we deduce that $\|d_0(\tilde{f_1}\circ f_1,r_{\alpha+\tilde{\alpha}})\|_{L^2(\Omega)}\leq 6|\lambda|^{\frac{1}{2}}$, and so $\|d_0(f_1\circ\tilde{f_1},\tilde{f_1}\circ f_1)\|_{L^2(\Omega)}\leq 12 |\lambda|^{\frac{1}{2}}$, and finally by mean value inequality $\|d_0(f\circ \tilde{f},\tilde{f}\circ f)\|_{L^2(\Omega)}\leq 48|\lambda|^{\frac{1}{2}}$.

\begin{rem} One could expect a converse inequality by using Moser's ideas \cite{Moser} to obtain a diffeomorphohism $h$ such that $\|d_{K-r}(hfh^{-1},r_\alpha)\|_{L^2(\Omega)}\ll \|d_{K}(f\circ \tilde{f},\tilde{f}\circ f)||_{L^2(\Omega)}$ and then deduce from Theorem \ref{lyapu} that $|\lambda|^{\frac{1}{2}}\ll \|d_{K}(f\circ \tilde{f},\tilde{f}\circ f)||_{L^2(\Omega)}$ for some $K$.
\end{rem}

The proof of Theorem \ref{principal} follows a ``KAM scheme'': in the same way as Arnold linearization Theorem \cite{Arnold} for a single diffeomorphism or Moser linearization theorem \cite{Moser} for commuting diffeomorphisms, we linearize the equation $hfh^{-1}=r_\alpha$ at $h=Id$, $f=r_\alpha$ so that a solution of the linear equation gives an approximate solution of the initial equation and thus define a conjugation $h$ such that $hfh^{-1}$ is closer to rotations than $f$. We prove that this can be achieved if the obstruction $\lambda$ is small enough by using the estimate given by Theorem \ref{lyapu}. Then we reiterate the process in order to conjugate $f$ to random diffeomorphisms $f_n$ closer and closer to rotations. The diophantine condition allows to control $C^k$ norms of the conjugations (up to some loss of derivatives phenomenomen, known problem classical to solve in these kind of KAM scheme), and the rotation number condition $\rho(f)=\alpha$ ensures that the diophantine condition is satisfied at each step of the process. Finally, if $\lambda=0$ we check that the sequence of conjugations converges and gives a conjugation between $f$ and $r_\alpha$, and if $\lambda\not=0$, we stop the process when $\lambda$ becomes large in front of $\mbox{dist}(f_n,r_\alpha)$ and it gives the wanted conjugation.\\

This scheme of the proof is smilar to the one in the paper of Dolgopyat and Krikorian \cite{Dolgopyat} where they prove an analog result on the sphere $S^d$ for $d\geq 2$ (though only the case $\lambda=0$).

\subsection{Lyapunov exponent of random product of matrices}

Our technics also apply to estimate the Lyapunov exponent of the product of i.i.d. random matrices $2\times 2$ close to rotation matrices, by studying the action on the projective line, identified to $\T$. And in this case we do not require a diophantine condition on the angle of the rotation but only a weak non degenerescence condition.\\

Let $||\cdot||$ be a norm in $\mathcal{M}_2(\R)$. Let $M$ be a random variable in $GL_2(\R)$. such that $\E[|\ln_+\|M\|]<+\infty$. It is a well known result of Kesten-Furstenberg \cite{Furstenberg-Kesten} that if $(M_n)_{n\in\N}$ is a sequence of independant copies of $M$, then the limit
$$\Lambda=\lim_{n\to\infty}\frac{\ln{\|M_{n-1}\cdots M_0\|}}{n}$$
exists almost sureley and does not depend on the alea. We call this number Lyapunov exponent of $M$.

For $\alpha\in\T$, we denote by $R_\alpha$ the rotation matrix of angulus $\pi\alpha$, that is to say $R_\alpha=\begin{pmatrix}
\cos \pi\alpha & -\sin \pi\alpha \\ \sin \pi\alpha &\cos \pi\alpha \end{pmatrix}$.

The following theorem is the analog of Theorem \ref{lyapu} for random product of matrices.

\begin{thm1}\label{mainproj}
	Let $\alpha$ be a random variable in $\T$ which does not belong almost surely to $\{0,\frac{1}{2}\}$. Let $M$ be a random variable in $SL_2(\R)$ of the form $M=R_\alpha+E$. Let $\ep=\E[||E||^3]^{\frac{1}{3}}$, that we assume to be finite, and let $\Lambda$ be the Lyapunov exponent of $M$. Then
	$$\Lambda=\frac{1}{8}\E\pa{\left|Z e^{i\pi\alpha}-\E[Ze^{i\pi\alpha}]\pa{\frac{1-e^{2i\pi \alpha}}{1-\E[e^{2i\pi \alpha}]}}\right|^2}+O(\ep^3)$$
	where $$Z=(a+d)+i(b-c)=\mbox{Tr}(E)+i \mbox{Tr}(ER_{\frac{1}{2}})$$
	(in particular, $\Lambda=O(\ep^2)$). If $\alpha$ is constant (i.e. non random), the formula simplifies itself and becomes
	$$\Lambda=\frac{1}{8}\E\cro{\left|Z-\E[Z]\right|^2}+O(\ep^3)=\frac{Var(Z)}{8}+O(\ep^3).$$
	The term $O(\ep^3)$ represents here a quantity bounded by $C\ep^3$ where $C$ is a constant depending only on $\alpha$ (and is actually uniformly bounded on the sets $\{\|d(\alpha,\{0,\frac{1}{2}\}\|_{L^2(\Omega)}\geq \mbox{const.}\}$)
\end{thm1}

\begin{rem}~
	
\begin{itemize}
	\item  In the general case $M\in GL_2(\R)$ (instead of $SL_2(\R)$), we can also obtain a Taylor expansion of its Lyapunov exponent $\Lambda$ by applying the Theorem to estimate the Lyapunov exponent $\widetilde{\Lambda}$ of $\widetilde{M}=M/\sqrt{\det(M)}$, since then $\Lambda=\widetilde{\Lambda}+\frac{1}{2}\E[\ln(\det(M))]$.
	\item As in Theorem \ref{lyapu}, the method can be generalized to obtain a Taylor expansion at any order, but it requires more restrictions on $\alpha$: to obtain an expansion at order $q$, $\alpha$ must not belong a.s. to $\{0,\frac{1}{q},\ldots,\frac{q-1}{q}\}$.
\item We can obtain from the theorem  an estimate of Figotin and Pastur \cite{Figotin} for the Lypunov exponent of a Schrodinger matrix with small random potential: if $M=\begin{pmatrix}E-gV&-1\\1&0\end{pmatrix}$, with $E=2\cos(\theta)\in]-2,2[-\{0\}$ and $V$ a random real variable having a third moment, then $M$ is conjugated to $R_\theta+gV\begin{pmatrix}1&\cot\theta\\0&0\end{pmatrix}$ and then by Theorem \ref{lyapu}, when $g$ tends to $0$ :
$$\Lambda=\frac{Var(V)}{8\sin^2\theta}g^2+O(g^3)=\frac{Var(V)}{2(4-E^2)}g^2+O(g^3).$$
\end{itemize}
\end{rem}

The following theorem is the analog of Theorem \ref{principal} for random product of matrices.

\begin{thm1}\label{mainproj2}
	Let $\mathcal{R}$ be the set of rotation matrices. For any $\delta>0$, there exists a neighborhood $\mathcal{U}$ of $\mathcal{R}$ in $SL_2(\R)$ such that for any random variable $M$ in $\mathcal{U}$ satisfying $||Tr(M)||_{L^2(\Omega)}\leq 2-\delta$, there exists $P\in SL_2(\R)$ such that
	$$\|d(PMP^{-1},\mathcal{R})\|_{L^2(\Omega)}\leq C\Lambda^{\frac{1}{2}},$$
where $\Lambda$ is the Lyapunov exponent of $M$ and $C$ is a constant depending only on the chosen norm on $\mathcal{M}_2(\R)$. Moreover, $\|P-I_2\|\leq C' \|d(M,\mathcal{R})\|_{L^2(\Omega)}$ for some $C'$ depending on $\delta$ and the norm.
\end{thm1}
From the proof it should not be difficult to explicit a constant $C$ for a given norm. The assumption $||Tr(M)||_{L^2(\Omega)}\leq 2-\delta$ gives a control of the ellipticity of $M$ in average, and should be seen as the analog of the the diophantine condition on $\rho(f)$ in the non linear case.

 We also deduce the same corollary as in the non linear case (with the same proof)

\begin{cor1}\label{commutateur}
	For any $\delta>0$, there exists a neighborhood $\mathcal{U}$ of $\mathcal{R}$ in $SL_2(\R)$ such that for any random variable $M$ in $\mathcal{U}$ satisfying $||Tr(M)||_{L^2(\Omega)}\leq 2-\delta$, if $\widetilde{M}$ is an independant copy of $M$ we have
	$$\E\cro{\|M\widetilde{M}-\widetilde{M}M\|^2}\leq C\Lambda,$$
	where $\Lambda$ is the Lyapunov exponent of $M$ and $C$ is a constant depending only on the chosen norm on $\mathcal{M}_2(\R)$.
\end{cor1}

From the proof it should not be difficult to obtain an explicit constant $C$ for a given norm. Moreover, by using compacity aguments in $\mathcal{M}_2(\R)$ we can deduce global results in more specific contexts, but then one can not hope to explicit the constants anymore without additional work. Here is an example of global result:
\begin{cor1}
Let $m$ be an integer and let $\delta$ and $C_0$ be two positive numbers, Then there exists $C>0$ such that for any matrices $A_1,\ldots,A_m$ in $SL_2(\R)$ satisfying $|Tr(A_i)|\leq 2-\delta$ (control of the ellipticity) and $\|A_i\|\leq C_0$ (control of the norm), we have  
	$$\sup_{i,j} \|A_i A_j -A_j A_i\|\leq C\Lambda^{\frac{1}{2}},$$
	where $\Lambda$ is the Lyapunov exponent of the uniformly distributed random matrix in $\{A_1,\ldots,A_m\}$.
\end{cor1}
\begin{proof}

	Let us consider $\Lambda$ as a function of $A_1,\ldots, A_m$ on $SL_2(\R)^m$. It is known by \cite{Bocker} that this function is continuous. In particular it is continuous on the compact subset 
	$$\mathcal{K}=\{(A_1,\ldots,A_m), \|A_i\|\leq C_0, |Tr(A_i)|\leq 2-\delta\}$$
	(the continuity of $\Lambda$ is actually a lot easier to prove on this subset $\mathcal{K}$ thanks to the ellipticity condition $|Tr(A_i)|\leq 2-\delta$).
	
	Moreover, if the function $\Lambda$ vanishes at a point $(A_1,\ldots,A_m)$ then by the classical Furstenberg Theorem \cite{Furstenberg} (and the ellipticity condition) the matrices $A_i$ commute. Thus there exists $P$ in $SL_2(\R)$ such that $PA_iP^{-1}$ is a rotation for every $i$, and using that $\|A_i\|\leq C_0$ and $|Tr(A_i)|\leq 2-\delta$  one can actually choose $P$ with a controled norm $\|P\|\leq C_1$ for some constant $C_1$ depending only on $C_0$ and $\delta$ (we leave this detail to the reader).
	
	Let $\mathcal{U}$ be the open set given by Corollary \ref{commutateur}, and let $$\mathcal{V}=\bigcup_{||P||\leq C_1}\left(P\mathcal{U}P^{-1}\right)^m\subset SL_2(\R)^m.$$ Then, $\Lambda$ is continuous and does not vanish on the compact set $\mathcal{K}\setminus \mathcal{V}$, hence $\Lambda\geq m$ for some $m>0$. Then:
	\begin{itemize}
		\item if $(A_1,...,A_m)\in \mathcal{V}$, there is $P$ in $Sl_2(\R)$ with $\|P\|\leq C_1$ such that $B_i=PA_iP^{-1}\in \mathcal{U}$ for every $i$, by Corollary \ref{commutateur} $\|B_i B_j -B_j B_i\|\leq C\Lambda^{\frac{1}{2}}$ for some constant $C$, and then $\|A_i A_j -A_j A_i\|\leq C'\Lambda^{\frac{1}{2}}$ for some new constant $C'=CC_1^2$
		\item if $(A_1,...,A_m)\notin \mathcal{V}$, then $\Lambda\geq m$ so $\|A_i A_j -A_j A_i\|\leq 2C_0^2\leq C\Lambda^{\frac{1}{2}}$ with $C=\frac{2C_0^2}{m^{\frac{1}{2}}}$.
	\end{itemize}
	
\end{proof}

\begin{rem}
	In the corollary above, one can actually obtain also a converse inequality $\sup_{i,j} \|A_i A_j -A_j A_i\|\geq c\Lambda^{\frac{1}{2}}$, by using that we can find $P$ with controlled norm and rotations matrices $R_i$ so that $\sup_i \|PA_iP^{-1}-R_i\|\ll \sup_{i,j} \|A_i A_j -A_j A_i\|$ and then by using Theorem \ref{mainproj} to get $\Lambda \ll \left(\sup_i \|PA_iP^{-1}-R_i\|\right)^2$. 
\end{rem}
\section{Preliminaries}\label{prelim}
\subsection{Some $C^k$ estimates}
We begin by state various estimates in $\mbox{Diff}_+^k(\T)$. All of them are classical estimates of KAM theory. Nevertheless, we give proofs in an appendix (section \ref{appen}).\\

 A key tool is the so called \textit{Kolmogorov inequality}.
\begin{prop1}(Kolmogorov inequality)\label{Kolmo}\\
 For any integers $j\leq k$ and for any $\varphi$ in $C^k(\T)$,
\begin{equation}\|\varphi\|_j\leq C\|\varphi\|_k^{j/k}\|\varphi\|_0^{1-j/k}.\end{equation}
where $C$ is a constant depending only on $k$.
\end{prop1}
The three following propositions give $C^k$ estimates of $gfg^{-1}$ when $f$ is a diffeomorphism close to a rotation $r_\alpha$ and $g$ is a diffeomorphism close to $Id$.
The first estimate allows to control the large $C^k$ norms of such a conjugation:
\begin{prop1}\label{estiK}
Let $f$, $g$ be in $\mbox{Diff}_+^k(\T)$ and let $\alpha$ be in $\T$ with $d_1(f,r_\alpha)\leq 1$ and $d_1(g,Id)\leq \frac{1}{2}$. Then :
$$d_k(gfg^{-1},r_\alpha)\leq C(d_k(f,r_\alpha)+d_k(g,Id)).$$
where $C$ is a constant depending only on $k$.
\end{prop1}
The assumption of the bound $1$ for $d_1(f,Id)$ is arbitrary and could be replace by any other number. In the same way the bound $\frac{1}{2}$ for $d_1(g,Id)$ could be replaced by any number less than $1$.\\

The second estimate bounds the distance between two conjugations in function of the distance between the cojugacies.
\begin{prop1}\label{esti0}
Let $f$, $g$ and $\tilde{g}$ be in $\mbox{Diff}_+^1(\T)$ and let $\alpha$ be in $\T$, with $d_1(f,r_\alpha)\leq 1$, $d_1(g,Id)\leq \frac{1}{2}$ and $d_1(\tilde{g},Id)\leq \frac{1}{2}$. Then:
$$ d_0(gfg^{-1},\tilde{g}f\tilde{g}^{-1})\leq C_0d_0(g,\tilde{g})$$
where $C_0$ is an absolute constant.
\end{prop1}
\begin{rem}
It is actually more generally possible to bound $d_k(gfg^{-1},\tilde{g}f\tilde{g}^{-1})$ in function of $d_k(g,\tilde{g})$, but we will not need it.
\end{rem}
The third estimate gives a classical linear approximation of $gfg^{-1}$
\begin{prop1}\label{esti1}
Let $k\geq 2$, let $f$, $g$ be in $\mbox{Diff}_+^2(\T)$ and let $\alpha$ be in $\T$. Writing $f=r_\alpha+\zetaup$, $g=Id+\etaup$ and denoting $\ep=\max(\|\zetaup\|_2,\|\etaup\|_2)$, we have
$$gfg^{-1}=r_\alpha+\left(\zetaup+\eta\circ r_\alpha-\eta\right)+R$$
where $R$ is a quadratic remainder satisfying $\|R\|_{1}\leq C\ep^2$ for some absolute constant $C$.
\end{prop1}
\begin{rem}
The $\ep^2$ upper bound can actually be replaced by the more precise term $\max(\|\zetaup\|_2,\|\etaup\|_2)\cdot\max(\|\zetaup\|_0,\|\etaup\|_0)$. There also exists a $C^k$ version of this estimate.
\end{rem}
We conclude with a last required estimate.
\begin{prop1}\label{estival}
Let $f$, $g$, $h$ be in $\mbox{Diff}_+^k(\T)$ with $d_k(h,Id)\leq 1$. Then :
$$d_k(f\circ h, g\circ h)\leq Cd_k(f,g).$$
where $C$ is a constant depending only on $k$.	
\end{prop1}
\begin{rem}
	Note that at the contrary of the previous propositions, we need to bound a large norm $d_k(h,Id)$, this is a strong assumption. Under the weak assumption $d_1(h,Id)\leq 1$ we actually have $d_k(f\circ h, g\circ h)\leq C(1+d_k(h,Id))d_k(f,g).$
\end{rem}

\subsection{Cohomological equation} 

We fix $r_\alpha=Id+\alpha$ a random rotation and $f=r_\alpha+\zetaup$ a perturbation of $r_\alpha$. We assume that $\alpha$ is a $(A,\sigma)$-diophantine. We will assume that $\sigma$ is an integer, in order to avoid the use of $C^k$-norms with $k$ non integer. It is obviously not a restriction since we can replace $\sigma$ by $[\sigma]+1$.\\
 
 We denote respectively by $T_0$ and $T$ the transfer operators of $r_\alpha$ and $f$. That is, for any map $\varphi:\T\rightarrow\R$,
 $$T_0\varphi=\E[\varphi\circ r_\alpha], \mbox{       } T\varphi=\E[\varphi\circ f].$$
Since $f$ is a perturbation of $r_\alpha$, $T$ is a perturbation of $T_0$. Note also that a measure $\mu$ is stationary for $f$ if and only if $\int\varphi d\mu=\int T\varphi d\mu$ for any map $\varphi\in C(\T)$.\\

The understanding of stationary measures is naturally related to the understanding of the cohomological equation $\varphi-T\varphi=\psi$.  Our main ingredient in our proofs is that the approximated cohomological equation  $\varphi-T_0\varphi=\psi$
 is easily solvable in $\varphi$ by Fourier methods, in the same way as in the classical deterministic case: the equation can be rewritten
$$\forall q\in \Z, \hat{\varphi}(q)(1-\E[e^{2i\pi q\alpha}])=\hat{\psi}(q).$$
For $q=0$ we get the obvious restriction $\hat{\psi}(0)=\int_\T\psi(x)dx=0$, and for $q\not=0$, if $q\alpha$ is not almost surely an integer (which is the case for $\alpha$ diophantine), then $\E[e^{2i\pi q\alpha}]\not=1$ and we obtain $\hat{\varphi}(q)=\frac{\hat{\psi}(q)}{1-\E[e^{2i\pi q\alpha}]}$. It leads us to define the following operator $U$: for $\psi:\T\rightarrow \R$,
$$U\psi(x)=\sum_{q\in\Z^*}\frac{\hat{\psi}(q)}{1-\E[e^{2i\pi q\alpha}]}e^{2i\pi qx}.$$
This operator apriori well defined at least for $\psi$ trigonometrical polynomial, gives the unique solution $\varphi$ if it exists to the equation $$\varphi-T_0\varphi=\psi-\int_\T\psi(x)dx$$
such that $\int_{\T}\varphi dx=0$.\\
 It is also convenient to define its adjoint $\overline{U}$ by
$$\overline{U}\psi(x)=\sum_{q\in\Z^*}\frac{\hat{\psi}(q)}{1-\E[e^{-2i\pi q\alpha}]}e^{2i\pi qx},$$
so that for any map trigonometric polynomials $\psi_1$ and $\psi_2$ we have 
$$\int_\T U\psi_1(x) \psi_2(x)dx=\int_\T \psi_1(x)\overline{U}\psi_2(x)dx.$$
The following lemma states that under the diophantine condition, $U$ and $\overline{U}$ are acutally well defined on sufficiently smooth maps, and are bounded up to some loss of derivative.

\begin{lem1}\label{cohomo}
Let $k_0=2\sigma+2$. Then the operators $U$ and $\overline{U}$ are well defined on $C^{k_0}(\T)$, and for any integer $k$, if $\psi\in C^{k+k_0}(\T)$ then $U\psi\in C^k(\T)$ and  $\displaystyle\|U\psi\|_{k}\leq \frac{1}{A^2}\|\psi\|_{k+k_0}$. The same estimate holds if we replace $U$ by $\overline{U}$.
\end{lem1} 

\begin{proof}
 It si enough to prove that for any integer $k$ the inequality  $\displaystyle\|U\psi\|_{k}\leq \frac{1}{A^2}\|\psi\|_{k+k_0}$ holds for any trigonometric polynomial $\psi$ (the same estimate for $\overline{U}$ follows by replacing $\alpha$ with $-\alpha$). To estimate $\|U\psi\|_{k}$ we are going to bound for $q\not=0$ the Fourier coefficient
$$|\widehat{U\psi}(q)|=\left|\displaystyle\frac{\hat{\psi}(q)}{1-\E[e^{2i\pi q\alpha}]}\right|.$$ The numerator can be bounded by above by 
\begin{equation}\label{num}
|\hat{\psi}(q)|\leq\frac{\|\psi\|_{k+k_0}}{(2\pi|q|)^{k+k_0}}.\end{equation}
To bound by below the denominator, we use that for any real number $x$, writing $x=k+\theta$ with $k\in\Z$ and $|\theta|=d(x,\Z)\leq \frac{1}{2}$ we have
$$1-\cos(2\pi x)=2\left(\sin(\pi x)\right)^2=2\left(\sin(\pi \theta)\right)^2\geq 2\left(\frac{2}{\pi}\pi\theta\right)^2=8 d(x,\Z)^2\geq d(x,\Z)^2,$$
hence by using the diophantine condition (\ref{dio}), 
\begin{equation}\label{den}\begin{disarray}{ll}
|1-\E[e^{2i\pi q\alpha}]|&\geq 1-\E[\cos(2\pi q\alpha)]\\
&\geq \E\cro{d(q\alpha,\Z)^2}\\
&\geq \frac{A^2}{|q|^{2\sigma}}.
\end{disarray}\end{equation}
Thus (\ref{num}) and (\ref{den}) give, using that $k_0=2\sigma+2$:
$$\displaystyle |\widehat{U\psi}(q)|\leq \frac{\|\psi\|_{k+k_0}}{(2\pi)^{k+k_0}A^2|q|^{k+2}}.$$
In consequence,
$$\displaystyle\|U\psi\|_{k}\leq \sum_{q\in\Z^*} |2\pi q|^{k}|\widehat{U\psi}(q)|\leq \frac{1}{(2\pi)^{k_0}A^2}\pa{\sum_{q\in\Z^*}\frac{1}{|q|^2}}\|\psi\|_{k+k_0}\leq \frac{1}{(2\pi)^2A^2}\frac{\pi^2}{3}\|\psi\|_{k+k_0}\leq \frac{1}{A^2}\|\psi\|_{k+k_0}.$$
\end{proof}
\section{Proof of Theorem \ref{lyapu}}

We fix a random rotation $r_\alpha$ and a perturbation $f=r_\alpha+\zetaup$, and we assume that $\alpha$ is $(A,\sigma)$-diophantine. The operators $T_0$, $T$, $U$ and $\overline{U}$ are defined as in previous section. We are going to obtain a Taylor expansion for the stationary measures of $f$ and the associated Lyapunov exponents.
\subsection{Estimate of the stationary measures}
\begin{prop1}\label{estimu}
If $\mu$ is a stationary measure for $f$, then:
$$\int_\T\varphi d\mu=\int_{\T}\varphi dx+O(\ep\|\varphi\|_{k_1})=\int_\T\varphi dx+\int_\T(\overline{U}\bar{\zetaup})\varphi 'dx+O(\ep^2\|\varphi\|_{k_2})$$
where $k_1=2\sigma+3$, $k_2=4\sigma+6$, $\bar{\zetaup}=\E[\zetaup\circ r_{-\alpha}]$ and $\ep=\E\left[\|\zetaup\|_{k_1}^2\right]^{\frac{1}{2}}$.
\end{prop1}
(As before $O(M)$ is a notation for a quantity bounded by $CM$ where $C$ is a constant depending only on $A$ and $\sigma$)
\begin{proof}
To prove the first equality of the statement, we start from the Taylor formula at order $0$: $\varphi\circ f=\varphi\circ r_\alpha+O(\|\zetaup\|_0\|\varphi\|_1)$, and we take the expectation, so
$$T\varphi=T_0\varphi+O(\ep\|\varphi\|_1).$$
Then, we use the invariance of $\mu$ :
$$\int_\T(\varphi-T_0\varphi) d\mu=O(\ep\|\varphi\|_1).$$
For $\psi$ in $C^{2\sigma+3}(\T)$, we apply the previous formula to $\varphi=U\psi$ and we get, thanks to Lemma \ref{cohomo} with $k=1$:
\begin{equation}\label{order 0}\int_\T\psi d\mu=\int_\T\psi dx+O(\ep\|\psi\|_{2\sigma+3}).\end{equation}
That gives the first equality.
\\ \\
To prove the second equality of the statement, we use this time a Taylor formula at order $1$ :
$$T\varphi=T_0\varphi+\E[(\varphi '\circ r_\alpha)\zetaup]+O(\ep^2\|\varphi\|_2).$$
Using the invariance of $\mu$, the first estimate (\ref{order 0}) and the inequality $\|uv\|_k\leq 2^k\|u\|_k\|v\|_k$ (consequence of Leibnitz formula), we get:
$$\begin{disarray}{ll}\int_\T(\varphi-T_0\varphi) d\mu&=\int_\T \E[(\varphi '\circ r_\alpha)\zetaup]d\mu+O(\ep^2\|\varphi\|_2)\\
&=\int_\T \E[(\varphi '\circ r_\alpha)\zetaup]dx+O(\ep^2\|\varphi\|_2+\ep\|\E[(\varphi '\circ r_\alpha)\zetaup]\|_{2\sigma+3})\\
&=\int_\T\varphi '\bar{\zetaup}dx+O(\ep^2\|\varphi\|_{2\sigma+4})
\end{disarray}$$
As before, for $\psi$ in $C^{4\sigma+5}(\T)$ we take $\varphi=U\psi$ to get, thanks to Lemma \ref{cohomo} with $k=2\sigma+4$:
$$\begin{disarray}{ll}\int_\T\psi d\mu&=\int_\T\psi dx+\int_\T (U\psi)'\bar{\zetaup}dx+O(\ep^2\|U\psi\|_{2\sigma+4})\\
&=\int_\T\psi dx+\int_\T \psi '(\overline{U}\bar{\zetaup})dx+O(\ep^2\|\psi\|_{4\sigma+6})
\end{disarray}$$
\end{proof}
\begin{rem}
We got that $\mu$ can be approximated by the density $h_0=1$ with accuracy $\ep$, and by the density $h_1=1-\overline{U}\bar{\zetaup}'$ with accuracy $\ep^2$ (in some sense to precise: we omit here the detail of the $C^k$-norms involved). We can easily generalize the method to have higher accuracy. Once defined an approximation $h_{n-1}$ with accuracy $\ep^{n-1}$, we write $T\varphi=T_0\varphi+T_1\varphi+\cdots+T_{n-1}\varphi+O(\ep^n\|\varphi\|)$ where $T_k\varphi=\frac{1}{k!}\E[(\varphi^{(k)}\circ r_\alpha)\zetaup^k]$. By a computation similar to the one in the proof we get $\int (\varphi-T_0\varphi)d\mu=\sum_{k=1}^{n-1}\int_{\T}\varphi \overline{T_k}h_{n-k} dx+O(\ep^n\|\varphi\|)$ where $\overline{T_k}\varphi=\frac{(-1)^k}{k!}\E[(\varphi^{(k)}\zetaup^k)\circ r_\alpha^{-1}]$. Then we apply to $\varphi=U\psi$ and we obtain that the density  $h_n=1+\sum_{k=1}^{n-1}\overline{U}~\overline{T_k}h_{n-k}$ approximate $\mu$ with accuracy $\ep^n$.
\end{rem}
\subsection{Estimate of the Lyapunov exponents}
Thanks to Proposition \ref{estimu} we can estimate the Lyapunov exponents of $f$:
\begin{prop1}\label{lyapu2} Let $k_0=4\sigma+7$.
If $\mu$ is a stationary probability for $f$ and $\lambda$ is the associated Lyapunov exponent, then
$$\lambda=-\frac{1}{2}\E\int_{\T}\left( \zetaup '-(\overline{U}\bar{\zetaup})'\circ r_\alpha+(\overline{U}\bar{\zetaup})'\right)^2dx+O(\ep^3)$$
where $\bar{\zetaup}=\E[\zetaup\circ r_{-\alpha}]$ and $\ep=\E[\|\zetaup\|_{k_0}^3]^{\frac{1}{3}}$.
\end{prop1} 
This will conclude the proof of Theorem \ref{lyapu}, setting $\eta=\overline{U}\bar{\zetaup}$.
\begin{proof}
Let $\eta=\overline{U}\bar{\zetaup}$, $g=Id-\eta$, $\tilde{f}=gfg^{-1}$ ($\|\eta\|_1=O(\ep)$ so $g$ is invertible if $\ep$ is small enough), $\tilde{\zetaup}=\tilde{f}-r_{\alpha}$ and $\tilde{\mu}=g_*\mu$. If $\varphi$ is in $C^{4\sigma+5}(\T)$, then thanks to Proposition \ref{estimu}, writing $\varphi \circ g=\varphi-\varphi'\etaup+O(\ep^2)$, we have, keeping the notations $k_1=2\sigma+3$ and $k_2=4\sigma+6$:
$$\begin{disarray}{ll}\int_\T\varphi d\tilde{\mu}&=\int_\T\varphi\circ g d\mu\\
&=\int_\T\varphi d\mu-\int_\T\varphi '\eta d\mu+O(\ep^2\|\varphi\|_2)\\
&=\left(\int_{\T}\varphi dx+\int_{\T}\varphi '\eta dx\right)-\int_{\T}\varphi '\eta dx+O(\ep^2\|\varphi\|_{k_2}+\ep\|\eta\|_{k_1}\|\varphi'\|_{k_1})\\
&=\int_\T\varphi dx+O(\ep^2\|\varphi\|_{k_2}).
\end{disarray}$$
where we used Lemma \ref{cohomo} to get $\|\eta\|_{k_1}=O(\|\bar{\zetaup}\|_{k_1+2\sigma+2})=O(\ep)$. Thus $\tilde{\mu}$ is ``$\ep^2$-close'' to Lebesgue measure.\\

The Lyapunov exponent $\lambda$ of $f$ associated to $\mu$ is equal to the Lyapunov exponent of $\tilde{f}$ associated to $\tilde{\mu}$ (this invariance of Lyapunov exponent by conjugation follows by taking the expectation and integrating with respect to $\mu$ the equality $\ln((gfg^{-1})')\circ g=\ln f'+(\ln g'\circ f-\ln g')$ ). We use this fact and the previous computation to estimate $\lambda$. We also use that by Proposition \ref{estiK}, $||\tilde{\zetaup}||_{k}=O(||\zetaup||_{k}+\|\etaup\|_k)$,  and that by Proposition \ref{esti1}, $\tilde{\zetaup}'=\left(\zetaup'-\eta'\circ r_\alpha+\eta'\right)+R$ with $\E[R^2]^{1/2}=O(\ep^2)$ . Then:
$$\begin{disarray}{ll}\lambda &=\E\int_\T\ln (1+\tilde{\zetaup}')d\tilde{\mu}\\
&=\E\int_\T (\tilde{\zetaup}'-\tilde{\zetaup}'^2/2)d\tilde{\mu}+O(\ep^3)\\
&=\E\int_{\T}(\tilde{\zetaup}'-\tilde{\zetaup}'^2/2)dx+O(\ep^3)\\ 
&=-\frac{1}{2}\E\int_{\T}\tilde{\zetaup}'^2 dx+O(\ep^3)\\
&=-\frac{1}{2}\E\int_{\T}( \zetaup '-\eta'\circ r_\alpha+\eta')^2dx+O(\ep^3).
\end{disarray}$$
\end{proof}
\begin{rem}
	We could avoid the conjugation by $g$ to estimate $\lambda$ and directly expand $\E\int\ln{ f '(x)}d\mu(x)$ using Proposition \ref{estimu}, but the method we have used has the advantage to make appear a main term clearly non-positive in the expansion of $\lambda$. Moreover, in the context of Theorem \ref{principal} this conjugation $g$ will correspond to the first step of the KAM scheme in order to conjugate $f$ to a diffeomorphism closer to rotations.
	\end{rem}

\section{Proof of Theorem \ref{principal}}
\subsection{Preliminaries}
We begin by introduce some convenient notations: if $u$ is a random variable valued in $C^k(\T)$, we set 
$$|||u|||_k=\E[\|u\|_k^2].$$

To avoid the profusion of constants, if $k$ is an integer we write $X\ll_k Y$ if $X\leq CY$ with $C$ a constant depending only on $A$, $\sigma$ and $k$, or simply $X\ll Y$ if $C$ depends only on $A$ and $\sigma$.\\

An other important tool is the smoothing operators, allowing to fix the loss of derivative phenomenom which will occur in the KAM scheme. Here we are going to simply use Fourier truncation, which does not give the optimal estimates but is sufficient for our purpose. So, for $\varphi:\T\rightarrow \R$ and $T\geq 0$ we denote
$$\left\{\begin{array}{l}\displaystyle
S_T\varphi(x)=\sum_{|p|\leq T}\hat{\varphi}(p)e^{2i\pi p x}\\ 
\displaystyle R_T\varphi(x)=\sum_{|p|>T}\hat{\varphi}(p)e^{2i\pi p x}.
\end{array}\right.$$
Then we have the standard Fourier estimates:
\begin{prop1}\label{Fourier}
For any integers $j$ and $k$ with $j<k$, we have
\begin{equation}\label{Kolmo1}\left\{\begin{array}{l}
\forall\varphi\in C^j(\T), \|S_T\varphi\|_k\ll_k T^{k-j+1}\|\varphi\|_j\\
\displaystyle\forall\varphi\in C^k(\T),\|R_T\varphi\|_j\ll_k \frac{\|\varphi\|_k}{T^{k-j-1}}.
\end{array}\right.\end{equation}
\end{prop1}
\subsection{First conjugation}
In this section we fix a random diffeomorphism $f=r_\alpha+\zetaup$ with $\alpha=\rho(f)$ diophantine of type $(A,\sigma)$, and $\lambda$ a Lyapunov exponent of $f$ associated to some stationary measure $\mu$. We assume that $f$ is valued in the open set
$$\mathcal{U}_0=\{h\in \mbox{Diff}_+^1(\T), |h'-1|<\frac{1}{2}\}.$$
In other words, $\mathcal{U}_0$ is the $\frac{1}{2}$-neighborhood of the set of rotations in $\mbox{Diff}_+^1(\T)$.
\begin{lem1}\label{conjugaison1}
Let $k_0=4\sigma+7$ and $r=2\sigma+2$. There exists $C_0>0$ depending only on $A$ and $\sigma$ so that $f$ is conjugated by a deterministic diffeomorphism $g=Id-\eta$ to $\tilde{f}=gfg^{-1}=r_\alpha+\tilde{\zetaup}$ such that either
$$|||\tilde{\zetaup}|||_0\leq 3|\lambda|^{\frac{1}{2}} \quad \mbox{or} \quad |||\tilde{\zetaup}|||_0\leq C_0|||\zetaup|||_{k_0}^{\frac{3}{2}},$$
with $\eta$ statisfying that for any integer $K\geq r$,
$$\|\etaup\|_{K-r}\ll_K |||\zetaup|||_{K}.$$
\end{lem1}
\begin{proof}
We begin with the same setting as in Proposition \ref{lyapu2}. First we set $\eta=\overline{U}\bar{\zetaup}$, which satisifes the inequality $\|\etaup\|_{K-r}\ll_K |||\zetaup|||_{K}$ by Lemma \ref{cohomo}. In particular $\|\etaup\|_{1}\ll |||\zetaup|||_{k_0}$ so we can assume $|||\zetaup|||_{k_0}$ small enough so that $\|\etaup\|_1 <\frac{1}{7}$ (if not, then $g=Id$ satisfies the conclusion of the statement). Then we set $g=Id-\eta$ which is invertible, $\tilde{f}=gfg^{-1}$,  $\tilde{\zetaup}=\tilde{f}-r_{\alpha}$ and $\tilde{\mu}=g_*\mu$.

Now, we follow the computation of the proof of Proposition \ref{lyapu2} with one slight difference: we cannot expand $\ln(1+\tilde{\zetaup}')$ at order $3$ because we do not have a good bound for the third moment of $||\zetaup||_1$. Instead we use that  for every $t$ in $]-1,1[$ we have $\ln(1+t)\leq t-\frac{1}{4}t^2$. We can apply this inequality to $t=\tilde{\zetaup} '$ because $f\in \mathcal{U}_0$ so $\tilde{f}'\leq \sup(f')\frac{\sup(g')}{\inf(g')}< (1+\frac{1}{2})\frac{1+\frac{1}{7}}{1-\frac{1}{7}}=2$ and so $-1<\tilde{\zetaup}'<1$. We get
$$\lambda=\E\int_\T \ln(1+\tilde{\zetaup}')d\tilde{\mu}\leq \E\int_\T (\tilde{\zetaup}'-\tilde{\zetaup}'^2/4)d\tilde{\mu}=-\frac{1}{4}\int_\T \tilde{\zetaup}'^2 dx +O(|||\zetaup|||_{k_0}^2)$$
hence there exists $C$ depending only on $A$ and $\sigma$ such that
$$\E\int_\T\tilde{\zetaup}'^2dx\leq 4|\lambda|+C|||\zetaup|||_{k_0}^3.$$
Next, we notice that for a fixed event, for every $a$, $b$, $|\tilde{\zetaup}(a)-\tilde{\zetaup}(b)|\leq \int_\T|\tilde{\zetaup}'|dx$, and since $\rho(\tilde{f})=\rho(f)=\alpha$, we have $\tilde{\zetaup}(b)=0$ for some $b$, and so $\|\tilde{\zetaup}\|_0\leq \int_\T|\tilde{\zetaup}'|dx$. Thus, by Cauchy-Schwarz, $\|\tilde{\zetaup}\|_0^2\leq\int_\T\tilde{\zetaup}'^2dx$, and taking the expectation we get 
$$|||\tilde{\zetaup}|||_0\leq \left(4|\lambda|+C|||\zetaup|||_{k_0}^{3}\right)^{\frac{1}{2}}\leq \left(\max(8|\lambda|,2C|||\zetaup|||_{k_0}^{3})\right)^{\frac{1}{2}}=\max\left(3|\lambda|^{\frac{1}{2}},\sqrt{2C}|||\zetaup|||_{k_0}^{\frac{3}{2}}\right),$$
which concludes the proof with $C_0=\sqrt{2C}$.
\end{proof}
In view of the dichotomy given by this lemma, we will say that ``$\lambda$ is an obstruction for the linearization of $f$ '' if $|\lambda|^{\frac{1}{2}}\geq \frac{C_0}{3}|||\zetaup|||_{k_0}^{\frac{3}{2}}$ where $C_0$ and $k_0$ are defined in the lemma. Thus, if $\lambda$ is an obstruction then one can find a conjugacy as stated in Theorem \ref{principal}, and if it is not an obstruction then $f$ is conjugated to a new random diffeomorphism $\tilde{f}$ closer to $r_\alpha$ and we can hope to iterate the process.
Though, we cannot use directly the lemma in an iterating process because of the loss of regularity in the inequality $|||\tilde{\zetaup}|||_0\leq C_0|||\zetaup|||_{k_0}^{\frac{3}{2}}$. We fix that by replacing the conjugation $g$ by a good $C^\infty$ approximation. In that way, there will be no loss of regularity anymore (at the cost of a less sharp bound). Precisely: 
\begin{lem1}\label{conjugaison2}
Let $k_0=4\sigma+7$ and $r=6\sigma+11$. If $\lambda$ is not an obstruction for $f$ then for any $T\geq 1$, $f$ is conjugated by a diffeomorphism $g_T=Id-\eta_T$ to $\tilde{f}_T=g_Tfg_T^{-1}=r_\alpha+\tilde{\zetaup}_T$ such that
$$\forall K\geq r, \left\{\begin{array}{l}\displaystyle|||\tilde{\zetaup}_T|||_{k_0}\ll_K 
T^{r}|||\zetaup|||_{k_0}^{\frac{3}{2}}+\frac{1}{T^{K-r}}|||\zetaup|||_K \\
\displaystyle |||\tilde{\zetaup}_T|||_K\ll_K T^r |||\zetaup|||_K\end{array}\right.$$
Moreover, 
$$
\forall K\geq r, \|\etaup_T\|_{K-r}\ll_K |||\zetaup|||_K
$$
\end{lem1}
\begin{proof}
Let $k_0=4\sigma+7$ and $s=2\sigma+2$. Let $g=Id-\etaup$ be the diffeomorphism given by Lemma \ref{conjugaison1}.  We set $\etaup_T=S_T\etaup$ and $g_T=Id-\etaup_T$. By Lemma \ref{conjugaison1} and Proposition \ref{Fourier} we have for $K\geq s+1$
\begin{equation}\label{state1}\|\etaup_T\|_{K-(s+1)}\ll_K \|\etaup\|_{K-s}\ll_K |||\zetaup|||_{K}.\end{equation}
Applying with $K=s+2\leq k_0$ we have $\|\etaup\|_1\ll |||\zetaup|||_{k_0}$, so we can assume $|||\zetaup|||_{k_0}$ small enough so that $\|\etaup\|_1\leq \frac{1}{2}$ (if not we set instead $g_T=Id$). Then $g_T$ is invertible and we can set $f_T=g_Tfg_T^{-1}=r_\alpha+\zetaup_T$. We also have for any $K\geq s+1$:
$$\|\etaup_T\|_{K}\ll_K T^{s+1}\|\etaup\|_{K-s}\ll_K T^{s+1}|||\zetaup|||_K,$$
so, by Proposition \ref{estiK}:
\begin{equation}\label{yop}|||\tilde{\zetaup}_T|||_K\ll_K|||\zetaup|||_{K}+\|\etaup_T\|_{K}\ll_K T^{s+1}|||\zetaup|||_K.\end{equation}
On another hand, since $\lambda$ is assumed not to be an obstruction for $f$ we have by Lemma \ref{conjugaison1}
$$|||gfg^{-1}-r_\alpha|||_0\ll |||\zetaup|||_{k_0}^{3/2},$$ 
and by Proposition \ref{esti0}, 
$$|||g_Tfg_T^{-1}-gfg^{-1}|||_0\ll\|g_T-g\|_0=\|R_T\etaup\|_0\ll_K \frac{1}{T^{K-s-1}}\|\etaup\|_{K-s}\ll_K \frac{1}{T^{K-s-1}}|||\zetaup|||_K.$$
\\
 The combination of the two last inequalities gives
 \begin{equation}\label{yap}
 |||\tilde{\zetaup}_T|||_{0}=|||g_Tfg_T^{-1}-r_\alpha|||_0\ll_K 
|||\zetaup|||_{k_0}^{\frac{3}{2}}+\frac{1}{T^{K-s-1}}|||\zetaup|||_K. \end{equation}
Finally, we write $\tilde{\zetaup}_T=S_T\tilde{\zetaup}_T+(\tilde{\zetaup}_T-S_T\tilde{\zetaup}_T)$ to use Proposition \ref{Fourier}, and then by using (\ref{yop}) and (\ref{yap}) we get
\begin{equation}\label{state2}|||\tilde{\zetaup}_T|||_{k_0}\ll_K T^{k_0+1}|||\tilde{\zetaup}_T|||_{0}+\frac{1}{T^{K-k_0-1}}|||\tilde{\zetaup}_T|||_{K}\ll_K T^{k_0+1}|||\zetaup|||_{k_0}^{\frac{3}{2}}+\frac{1}{T^{K-k_0-s-2}}|||\zetaup|||_K.\end{equation}
Thus, with $r=k_0+s+2=6\sigma+11$, (\ref{state1}), (\ref{yop}) and (\ref{state2}) give all the estimates claimed in the statement.
\end{proof}
\subsection{KAM iteration}
Now we begin the KAM scheme by iterating the conjugation process given by Lemma \ref{conjugaison2}. We fix $k_0$ and $r$ the numbers given by the Lemma \ref{conjugaison2}, and we fix a sequence of numbers $(T_n)_{n\in\N}$. We initialize the construction with $f_0=f$, $\zetaup_0=\zetaup$. Then, assuming that $f_{n-1}=r_\alpha+\zetaup_{n-1}$ is defined, if we have the two conditions

\begin{enumerate}\label{test}
	\item $f_{n-1}\in \mathcal{U}_0$ a.s.,
	\item $\lambda$ is not an obstruction for $f_{n-1}$, that is $|\lambda|^{\frac{1}{2}}\leq \frac{C_0}{3}|||\zetaup_{n-1}|||_{k_0}^{\frac{3}{2}}$,
\end{enumerate}
then Lemma \ref{conjugaison2} applies, so that by choosing $T=T_n$ we get a conjugation $g_{n-1}=Id-\etaup_{n-1}$ and a random diffeomorphism $f_{n}=g_{n-1}f_{n-1}g_{n-1}^{-1}=r_\alpha+\zetaup_{n}$ satisfying for $K\geq r$
$$\left\{\begin{disarray}{l}
|||\zetaup_{n}|||_K\ll_K T_n^{r}|||\zetaup_{n-1}|||_K\\
|||\zetaup_{n}|||_{k_0}\ll_K 
T_n^{r}|||\zetaup_{n-1}|||_{k_0}^{\frac{3}{2}}+\frac{1}{T_n^{K-r}}|||\zetaup_{n-1}|||_K
\end{disarray}\right.$$
and 
$$\|\eta_{n-1}\|_{K-r}\ll_K T_n^{r}|||\zetaup_{n-1}|||_K.$$
If one of the two conditions is not satisified, then we stop the process.
Thus we get a sequence of random diffeomorphisms $(f_n)_{n<N}$ where $N\in \N\cup\{+\infty\}$.

The choice of $T_n$ we do is the following: $T_n=2^{Q^n}$ where $Q$ is any number in $(1,\frac{3}{2})$. With this choice, we prove that the large $C^k$-norms of $\zetaup$ do not grow up too fast while the small $C^k$-norms decrease quickly. Note that in the sequel we consider $Q$ as fixed, for exemple $Q=\frac{4}{3}$, so we will not explicit the dependence of the constants in $Q$.
\begin{lem1}\label{induc}
	There exists integers $p$ and $K_0$ depending only on $\sigma$ such that for any $K\geq K_0$, if $\ep=|||\zetaup|||_{K}$ is small enough then for any $n<N$, 
	$$\left\{\begin{array}{l}
	|||\zetaup_n|||_K\ll_K T_{n}^p\ep\\
	\displaystyle|||\zetaup_n|||_{k_0}\ll_K \frac{1}{T_{n}^{K-p}}\ep
	\end{array}\right.$$
\end{lem1}
\begin{proof}
	There exists a constant $C$ depending only on $A$, $\sigma$ and $K$ such that for any $n<N$
	$$
	\left\{\begin{array}{l}
	\displaystyle |||\zetaup_{n}|||_{K}\leq C T_n^{r}|||\zetaup_{n-1}|||_{K}\\
	\displaystyle |||\zetaup_{n}|||_{k_0}\leq C\left( 
	T_n^{r}|||\zetaup_{n-1}|||_{k_0}^{\frac{3}{2}}+\frac{1}{T_n^{K-r}}|||\zetaup_{n-1}|||_{K}\right)
	\end{array}\right.
	$$
	By iteration of the first inequality we have for any $n\geq 1$:
	$$|||\zetaup_n|||_{K}\leq C^n(T_{n}\cdots T_1)^r|||\zetaup_0|||_{K}\leq C^{n} 2^{r(Q+Q^2+\ldots+Q^{n})}\ep\leq C^{n} 2^{\frac{rQ}{Q-1} Q^{n}}\ep,$$
	hence $|||\zetaup_n|||_K\ll_K T_n^s\ep$ where $s=\frac{2rQ}{Q-1}$. That proves the first part of the statement if $p\geq s$.\\
	
	Let $\ep_n=|||\zetaup_n|||_{k_0}$. Using in the second inequality that $|||\zetaup_{n-1}|||_K\ll_K T_n^{s}\ep$, we obtain, up to modifying the constant $C$:
	$$\ep_{n}\leq C \left(T_n^{r}\ep_{n-1}^{\frac{3}{2}}+\frac{1}{T_n^{K-p}}\ep\right).$$
	where we have set $p=r+s$. If $K$ is large enough and $\ep$ small enough, we are going to prove by induction that for every $n<N$, 
	\begin{equation}\label{recu}\ep_n\leq \frac{2C\ep}{T_{n}^{K-p}}.\end{equation}
	It holds for $n=0$ if $C\geq 2^K$, what we can assume up to changing $C$ one more time. Now, for $n<N$ let us assume that $\ep_{n-1}\leq \frac{2C\ep}{T_{n-1}^{K-p}}$. Then if $\ep$ is small enough we have
	$$\ep_{n-1}^{\frac{3}{2}}\leq \frac{1}{T_{n-1}^{\frac{3}{2}(K-p)}}(2C\ep)^{\frac{3}{2}}\leq \frac{1}{T_n^{\frac{3}{2Q}(K-p)}}\ep,$$
	and so
	$$\ep_{n}\leq C\ep\left(\frac{1}{T_n^{\frac{3}{2Q}(K-p)-r}}+\frac{1}{T_n^{K-p}}\right),$$
	which implies that 
	$$\ep_{n}\leq\frac{2C\ep}{T_n^{K-p}}$$
	provided that $\frac{3}{2Q}(K-p)-r\geq K-p$, or equivalently (since $\frac{3}{2Q}>1$)
	$$K\geq p+s\frac{1}{\frac{3}{2Q}-1}.$$
	If it is satisfied then (\ref{recu}) is proved by induction for any $n<N$. That concludes the proof of the lemma, choosing $K_0=\lceil p+s\frac{1}{\frac{3}{2Q}-1}\rceil $.
\end{proof}
In the sequel we fix the integer $K_0$ given by Lemma \ref{induc}, and an integer $K\geq K_0$.
\begin{lem1}\label{upgrade}
	There exists $q$ depending only on $\sigma$ such that if $\ep=|||\zetaup|||_K$ is small enough then
	for any $n<N$, $|||\zetaup_n|||_{K-q}\ll_K \frac{1}{T_{n}}\ep$ and $||\etaup_n||_{K-q}\ll_K \frac{1}{T_{n}}\ep$
\end{lem1}
\begin{proof}
	Let $p$ as in previous lemma and let $K\geq K_0$. If $\ep$ is small enough we have $
	|||\zetaup_n|||_K\ll_K T_{n}^p\ep$ and $|||\zetaup_n|||_{0}\ll_K \frac{1}{T_{n}^{K-p}}\ep$, so by Kolomogorov inequality (Proposition \ref{Kolmo}), for any $k\leq K$ we have
	$$|||\zetaup_n|||_{K-k}\ll_K|||\zetaup_n|||_0^{\frac{k}{K}}|||\zetaup_n|||_K^{\frac{K-k}{K}}\ll_K \frac{\ep}{T_{n}^\tau}$$
	with
	$$\tau=\frac{k}{K}(K-p)-\left(\frac{K-k}{K}\right)p=k-p$$
	In particular, $|||\zetaup_n|||_{K-q}\ll_K \frac{1}{T_{n}}\ep$  if $q\geq p+1$, and $||\etaup_n||_{K-q}\ll_K |||\zetaup_n|||_{K-q+r}\ll_K \frac{1}{T_{n}}\ep$ if $q-r\geq p+1$. So we get the result with $q=p+1+r$.
	
\end{proof}

Now we consider the compositions $h_n=g_{n-1}\circ\cdots\circ g_0$, so that $f_n=h_{n}fh_{n}^{-1}$. The diffeomorphisms $h_n$ satisfy the following estimates:
\begin{lem1}\label{finish}
	Let $q$ as in previous lemma. If $\ep=|||\zetaup|||_K$ is small enough then for any $n<N$ $d_{K-q}(h_n,Id)\ll_k \ep$ and $\sum_{n<N} d_{K-q}(h_n,h_{n-1})\ll_K \ep$.
\end{lem1}

\begin{proof}
	
	Let $\delta_n=d_{K-q}(h_n,Id)$. For a fixed $n$ let us assume that $\delta_j\leq 1$ for $j=0,..,n-1$. Then, by Proposition \ref{estival} and Lemma \ref{upgrade},
	$$d_{K-q}(h_n,h_{n-1})\ll_K d_{K-q}(g_n,Id)\ll_K \frac{\ep}{T_{n}},$$
	and so 
	$$\delta_n\leq \sum_{j<n}d_{K-q}(h_j,h_{j-1})\ll_K \ep.$$
	So if $\ep$ is small enough we get $\delta_n\leq 1$. Thus we get by induction that $\forall n<N, \delta_n\leq 1$. In particular the estimates above hold for every $n$, and the result follows.
\end{proof}

We are now ready to finish the proof of Theorem \ref{principal}.

\begin{proof}(Theorem \ref{principal})\\
	We fix $K_0$ and $q$ as above, an integer $K\geq K_0$, we assume that $\ep=|||\zetaup|||_K$ is small enough so that the lemmas above apply, and we also assume that $|f'-1|\leq \frac{1}{4}$. We separate the cases $N=+\infty$ and $N<+\infty$.\\
	\begin{itemize}
		\item  If $N=+\infty$, then $\sum_n d_{K-q}(h_n,h_{n-1})\ll_K \ep$ hence
		$(h_n)_{n\in\N}$ converges in $\mbox{Diff}_+^{K-q}(\T)$ to a limit $h$
		satisfying $d_{K-q}(h,Id)\ll_K \ep$. In particular if $\ep$ is small enough $h$ is invertible and $hfh^{-1}=\lim_n h_nfh_n^{-1}=\lim_n f_n=r_\alpha$ almost surely.\\
		
		\item If $N<+\infty$. Then, $f_{N-1}=h_{N-1}fh_{N-1}^{-1}$ with $d_{K-q}(h_{N-1},Id)\ll_K \ep$. Morever, one of the two coditions stated at the beginning of the section does not hold for $f_{N-1}$, that is, either $f_{N-1}\notin \mathcal{U}_0$ or $\lambda$ is an obstruction for $f_{N-1}$. Since $|f'-1|\leq \frac{1}{4}$ and $|h_{N-1}'-1|\ll_K \ep$, we deduce that the  condition $f_{N-1}\in \mathcal{U}_0$ is satisfied if $\ep$ is small enough. So it means that   $\lambda$ is an obstruction for $f_{N-1}$, that is $|\lambda|^{\frac{1}{2}}\geq \frac{C_0}{3}\ep_n^{\frac{3}{2}}$. Then Lemma \ref{conjugaison1} gives a diffeomorphism $g$ satisfying $d_{K-q}(g,Id)\ll_K\ep$ conjugating $f_{N-1}$ to  $\tilde{f}=r_\alpha+\tilde{\zeta}$ such that $|||\tilde{\zetaup}|||_0\leq 3|\lambda|^{\frac{1}{2}}$, and then the conjugation $h=g\circ h_{N-1}$ satisfies the conclusion of Theorem \ref{principal}.
	\end{itemize}
	
	Choosing $\overline{\ep}$ in $(0,\frac{1}{2})$ so that the lemmas above and the final argument apply for $|||\zetaup|||_K\leq \overline{\ep}$, we get the conclusion of Theorem \ref{principal} for any random diffeomorphism $f$ such that $\rho(f)$ is $(A,\sigma)$-diophantine and valued in the open set
	$$\mathcal{U}=\left\{h\in \mbox{Diff}_+^K(\T), d_{K}(h,\mathcal{R})< \frac{\overline{\ep}}{2}\right\},$$
	where $\mathcal{R}$ is the set of rotations: for such a $f$, we obviously have $|f'-1|\leq \frac{1}{4}$, and $d_K(f,r_\beta)<\frac{\overline{\ep}}{2}$ for some $\beta$ so actually $|\beta-\alpha|<\frac{\overline{\ep}}{2}$ with $\alpha=\rho(f)$, so $d_K(f,r_\alpha)<\overline{\ep}$ and in particular $|||\zetaup|||_K\leq \overline{\ep}$. Hence the argument above applies to $f$ and gives the conjugation stated in Theorem \ref{principal}.
	
\end{proof}

\section{Random products of matrices (Theorems \ref{mainproj} and \ref{mainproj2})}
\subsection{Generalities}
We consider $\mathcal{M}_2(\R)$ equipped with any norm $||\cdot||$. By identifying the complex plane with $\R^2$, any matrix $M$ in $\mathcal{M}_2(\R)$ naturally acts on $\C$.\\

We denote by $\mathcal{T}$ the space of trigonometrical polynomials $p:\T\rightarrow\R$, generated by the maps $x\mapsto \cos(2k\pi x)$ and $x\mapsto \sin(2k\pi x)$. We denote by $\mathcal{T}_n$ the space of trigonometrical polynomials of $\mathcal{T}$ of degree at most $n$. We fix a norm $|| \cdot ||$ on $\mathcal{T}$.\\

To any $M$ in $GL_2(\R)$ we naturally associate a diffeomorphism $f_M$ of $\T$ by 
$$e^{i\pi f_M(x)}=\frac{M(e^{i\pi x})}{|M(e^{i\pi x})|}.$$

We admit the following elementary lemma.
\begin{lem1}\label{equiv}
	There exists a constant $A_0>0$ depending only of the norm on $\mathcal{M}_2(\R)$ such that for any $M$ in $SL_2(\R)$ and $\alpha$ in $\T$, 
	$$\frac{1}{A_0}d_0(f_M,r_\alpha)\leq ||M-R_\alpha||\leq A_0d_0(f_M,r_\alpha).$$
\end{lem1}

In particular, if $M$ is a perturbation of $R_\alpha$ of order $\ep$, then $f_M$ is a perturbation of $r_\alpha$ of order $\ep$. The next lemma specifies the form of the perturbation:
\begin{lem1}\label{proj}
	If $M=R_\alpha+E$ then, writing $f_M=r_\alpha+\zetaup$ we can write $\zetaup=\zetaup_1+\zetaup_2+\zetaup_3$ where $\zetaup_1\in \mathcal{T}_1$ and $||\zetaup_1||=O(\|E\|)$, $\zetaup_2\in \mathcal{T}_2$ and $||\zetaup_2||=O(\|E\|^2)$, $\zetaup_3\in C^\infty(\T)$ and $||\zetaup_3||_1=O(\|N\|^3)$. Moreover,
	$$\zetaup_1(x)=\frac{1}{\pi}\mbox{Im}\pa{E(e^{i\pi x})e^{-i\pi(x+\alpha)}}$$
\end{lem1}
\begin{proof}
	From $e^{i\pi f_M(x)}=\frac{M(e^{i\pi x})}{|M(e^{i\pi x})|}$ we obtain the formula
	$$\zetaup(x)=\frac{1}{i\pi}\ln\pa{\frac{1+E(e^{i\pi x})e^{-i\pi(x+\alpha)}}{|1+E(e^{i\pi x})e^{-i\pi(x+\alpha)}|}},$$
	where the (complex) logarithm is well defined for $||E||$ small. Then the result follows by doing Taylor expansions.
\end{proof}
The following lemma is a counterpart of the previous lemma when $\alpha=0$ that we will use to create a conjugation matrix in the proof of Theorem \ref{mainproj2}.
\begin{lem1}\label{projR}
	If $\zetaup$ belongs to $\mathcal{T}_1$, then one can find $M$ in $SL_2(\R)$ such that $||M-I_2||=O(||\zetaup\|)$ and
	$$f_M(x)=x+\zetaup(x)+O(\|\zetaup\|^2).$$
\end{lem1}	
\begin{proof}
	By assumption, $\zetaup(x)=A+B\cos(2\pi x)+C\sin(2\pi x)$ for some $A,B,C$. Let us set $M=I_2+E$ with $E=\begin{pmatrix}a & b \\ c & d\end{pmatrix}$, where $a,b,c$ have to be chosen, and $d$ is determined so that $\det M=1$. Since $\det(M)=1+ Tr(E)+O(||E||^2)$, in particular $d=-a+O(||E||^2)$. From Lemma \ref{proj} and a simple computation, we have
	$$\begin{disarray}{ll}f_M(x)&=x+\frac{1}{\pi}\mbox{Im}\pa{E(e^{i\pi x})e^{-i(\pi x+\alpha)}}+O(||E||^2)\\
	&=x+\frac{c-b}{\pi}+\frac{c+b}{\pi}\cos(2\pi x)+\frac{d-a}{\pi} \sin(2\pi x)+O(||E||^2)\\
	&=x+\frac{c-b}{\pi}+\frac{c+b}{\pi}\cos(2\pi x)-\frac{2a}{\pi} \sin(2\pi x)+O(||E||^2).\end{disarray}$$
	By chosing $a,b,c$ so that $c-b=\pi A$, $c+b=\pi B$ and $-2a=\pi C$, we obviously have $||E||=O(||\zetaup||)$ and so $f_M(x)=x+\zetaup(x)+O(||\zetaup||^2)$.
\end{proof}	

\begin{lem1}\label{Lyapuproj}
	Let $M$ be a random matrix in $SL_2(\R)$ with $\E[\ln_+ ||M||]<+\infty$, and let $\Lambda$ be the Lyapunov exponant of $M$. Then there exists a stationary measure $\mu$ of the random diffeomorphism $f_M$ so that the corresponding Lyapunov exponent $\lambda(\mu)$ satisfies $\Lambda=-\frac{1}{2}\lambda(\mu)$.
\end{lem1}
\begin{proof}
	Since $M\in SL_2(\R)$, we have for every $\theta$ and $\theta'$ in $\T$ $$\det(M(e^{i\pi\theta}),M(e^{i\pi\theta'}))=\det(e^{i\pi\theta},e^{i\pi\theta'}),$$ that we can rewrite 
	$$|\sin(\pi(\theta-\theta '))|=|M(e^{i\pi\theta})|~|M(e^{i\pi\theta '})|~|\sin(\pi(f_M(\theta)-f_M(\theta ')))|,$$ 
	which leads to 
	$$1=|M(e^{i\pi\theta})|^2~|f_M'(\theta)|.$$
	It is well known that there exists a stationary measure $\mu$ such that we have $\Lambda=\E\int_{\T}\ln |M(e^{i\pi\theta})|d\mu(\theta)$ (see for exemple \cite{Furstenberg}), so the result follows.
\end{proof}

\subsection{Proof of Theorem \ref{mainproj}}
We fix a random variable $\alpha$ in $\T$ and a random matrix $M=R_\alpha+E$ of $SL_2(\R)$. We naturally get a random diffeomorphism $f_M=r_\alpha+\zetaup$ of $\T$, and Lemma \ref{proj} gives a decomposition $\zetaup=\zetaup_1+\zetaup_2+\zetaup_3$.

We assume that $\alpha$ does not belong almost surely to $\{0,\frac{1}{2}\}$. So $||d(2\alpha,\Z)||_{L^2(\Omega)}\geq \delta$ for some $\delta>0$. In the sequel a term $O(M)$ means a term bounded by $CM$ with $C$ depending only on $\delta$ (and the chosen norms on $\mathcal{T}$ and $\mathcal{M}_2(\R)$).

We keep the notations of the previous sections for the operators $T, T_0, U$ and $\overline{U}$, that is to say $T\varphi(x)=\E[\varphi\circ f_M(x)]$, $T_0\varphi(x)=\E[\varphi\circ r_\alpha(x)]$, $U\varphi(x)=\sum_{q\in\Z^*}\frac{\hat{\varphi}(q)}{1-\E[e^{2i\pi q\alpha}]}e^{2i\pi qx}$
and $\overline{U}\varphi(x)=\sum_{q\in\Z^*}\frac{\hat{\varphi}(q)}{1-\E[e^{-2i\pi q\alpha}]}e^{2i\pi qx}$.

\begin{lem1}\label{cohomoproj}
	The operators $U$ and $\overline{U}$ are well defined and bounded on $\mathcal{T}_2$. Moreover, $||U||$ and $||\overline{U}||$ can be bounded by a constant depending only on $\delta$ (and the norm $||\cdot||$ on $\mathcal{T}_2$).
\end{lem1}
\begin{proof}
	The operators $U$ and $\overline{U}$ are well defined on $\mathcal{T}_2$ since the denominators $1-\E[e^{2i\pi q\alpha}]$ do not vanish for $q=-2,-1,1,2$ thanks to the assumption that  $\alpha$ does not belong almost surely to $\{0,\frac{1}{2}\}$. Theses operators are automatically bounded since $\mathcal{T}_2$ is finite dimensional. Finally the uniform bound of $||U||$ and $||\overline{U}||$ follows from the inequality $|1-\E[e^{2i\pi q\alpha}]|\geq 8\E\cro{d(q\alpha,\Z)^2}$ (obtained in the proof of Lemma \ref{cohomo}) applied to $q=-2,-1,1,2$.
\end{proof}
\begin{lem1}\label{LyapuProj}
	$$\Lambda=\frac{1}{4}\E\int_{\T} \pa{\zetaup_1 '+(\overline{U}\bar{\zetaup}_1)'-(\overline{U}\bar{\zetaup}_1)'\circ r_\alpha}^2dx+O(\ep^3).$$
	where $\bar{\zetaup}_1=\E[\zetaup_1\circ r_\alpha^{-1}]$ with $\zetaup_1$ given by Lemma \ref{proj}, and $\ep=\E[||E||^3]^{\frac{1}{3}}$.
\end{lem1}
\begin{proof}
	
	By Lemma \ref{Lyapuproj}, we have $\Lambda=-\frac{1}{2}\lambda(\mu)$ for some stationary probability measure $\mu$ on $\T$. If $\alpha$ is diophantine, the expansion in the statement is a consequence of Proposition \ref{lyapu}. We are going to check that the estimate is still valid without diophantine assumption by mimicking the proof of Proposition \ref{lyapu}, noticing that we only need to estimate $\mu$ on trigonometrical polynomials of small degrees, and so we only need the boundedness of $U$ on $\mathcal{T}_2$ given by Lemma \ref{cohomoproj}.
	\begin{itemize}
		\item For every $\psi$ in $\mathcal{T}_2$, with $\varphi=U\psi$ ($\in \mathcal{T}_2$) we have
		$$\int_{\T}\psi d\mu-\int_{\T}\psi dx=\int_{\T}(\varphi-T_0\varphi) d\mu=\int_{\T}(T\varphi-T_0\varphi) d\mu=O(\ep\|\varphi\|)=O(\ep\|\psi\|)$$
		\item For every $\psi$ in $\mathcal{T}_1$, with $\varphi=U\psi$ ($\in \mathcal{T}_1$) we have
		$$\begin{disarray}{ll}\int_{\T}\psi d\mu-\int_{\T}\psi dx&=\int_{\T}(T\varphi-T_0\varphi) d\mu\\
		&=\int_{\T}\E[(\varphi '\circ r_\alpha)\zetaup]d\mu+O(\ep^2\|\varphi\|)\\
		&=\int_{\T}\E[(\varphi '\circ r_\alpha)\zetaup_1]dx+O(\ep^2\|\varphi\|)\\
		&=\int_{\T}\varphi'\bar{\zetaup}_1dx+O(\ep^2\|\varphi\|)\\
		&=\int_{\T}\psi '\overline{U}\bar{\zetaup}_1 dx+O(\ep^2\|\psi\|)\end{disarray}$$
		(for the third equality we used that $(\varphi '\circ r_\alpha)\zetaup_1$ belongs to $\mathcal{T}_2$)
		\item Denoting $\etaup=\overline{U}\bar{\zetaup}_1$ ($\in\mathcal{T}_1$), $g=Id-\etaup$ and $\tilde{\mu}=g_*\mu$, we have for $\psi$ in $\mathcal{T}_2$
		$$\int_{\T}\psi d\tilde{\mu}=\int_{\T}\psi d\mu+O(\ep\|\psi\|)=\int_{\T}\psi dx+O(\ep\|\psi\|)$$
		and for $\psi$ in $\mathcal{T}_1$,
		$$\int_{\T}\psi d\tilde{\mu}=\int_{\T}\psi d\mu-\int_{\T}\psi ' \overline{U}\bar{\zetaup}_1 d\mu+O(\ep^2\|\psi\|)=\int_{\T}\psi dx+O(\ep^2\|\psi\|)$$
		\item Denoting $\tilde{f}=g\circ f_M\circ g^{-1}=r_\alpha+\tilde{\zetaup}$ ($g$ invertible if $\ep$ is small enough since $||\etaup||=O(\ep)$),
		by using the decomposition $\zetaup=\zetaup_1+\zetaup_2+\zetaup_3$ and Taylor expansions we can write $\tilde{\zetaup}=\tilde{\zetaup}_1+\tilde{\zetaup}_2+\tilde{\zetaup}_3$ with
		$$\left\{\begin{disarray}{l}
		\tilde{\zetaup}_1=\zetaup_1-\etaup\circ r_{\alpha}+\etaup,~ \tilde{\zetaup}_1\in 
		\mathcal{T}_1,~||\tilde{\zetaup}_1||=O(\max(\|E\|,\|\etaup\|))\\
		\tilde{\zetaup}_2\in \mathcal{T}_2,~||\tilde{\zetaup}_2||=O(\max(\|E\|^2,\|\etaup\|^2))\\
		||\tilde{\zetaup}_3||_1=O(\max(\|E\|^3,\|\etaup\|^3))
		\end{disarray}\right..$$
		\item We conclude:
		$$\begin{disarray}{ll}\lambda(\mu)&=\E\int_{\T}\ln \tilde{f}'d\tilde{\mu}\\
		&=\E\int_{\T}\tilde{\zetaup}_1'd\tilde{\mu}+\E\int_{\T}\tilde{\zetaup}_2'd\tilde{\mu}-\frac{1}{2}\E\int_{\T}\tilde{\zetaup}_1'^2d\tilde{\mu}
		+O(\ep^3)\\
		&=-\frac{1}{2}\int_{\T}\tilde{\zetaup}_1'^2dx+O(\ep^3)
		\end{disarray},$$
		from which the result follows since $\Lambda=-\frac{1}{2}\lambda(\mu)$
	\end{itemize}
\end{proof}

We can deduce Theorem \ref{mainproj} by a serie of simple computations. Starting from the equality $E(e^{i\pi x})=\frac{1}{2}(Ze^{i\pi x}+Z' e^{-i\pi x})$
with $Z=(a+d)+i(c-b)$ and $Z'=(a-d)+i(b+c)$, we successively obtain (using Lemma \ref{proj})
\begin{itemize}
	\item $\displaystyle \zetaup_1(x)=\frac{1}{\pi}\mbox{Im}\pa{E(e^{i\pi x})e^{-i\pi(x+\alpha)}}=\frac{1}{2\pi}\mbox{Im}\pa{Z e^{i\pi(2x+\alpha)}}+\mbox{constant}$
	\item $\displaystyle \bar{\zetaup}_1(x)=\frac{1}{2\pi}\mbox{Im}\pa{\E[Z e^{-i\pi\alpha}] e^{2i\pi x}}+\mbox{constant}$
	\item $\displaystyle \overline{U}\bar{\zetaup}_1(x)=\frac{1}{2\pi}\mbox{Im}\pa{\frac{\E[Ze^{-i\pi\alpha}]}{1-\E[e^{-2i\pi \alpha}]}e^{2i\pi x}}$
	\item $\displaystyle \pa{\zetaup_1 +\overline{U}\bar{\zetaup}_1-\overline{U}\bar{\zetaup}_1\circ r_\alpha}(x)=\frac{1}{2\pi}\mbox{Im}\pa{X e^{2i\pi x}}+\mbox{constant}$\\
	where
	$\displaystyle X=Z e^{i\pi\alpha}+\frac{\E[Ze^{-i\pi\alpha}]}{1-\E[e^{-2i\pi \alpha}]}-\frac{\E[Ze^{-i\pi\alpha}]}{1-\E[e^{-2i\pi \alpha}]}e^{2i\pi\alpha}$,
	\item $\displaystyle \pa{\zetaup_1 '+(\overline{U}\bar{\zetaup}_1)'-(\overline{U}\bar{\zetaup}_1)'\circ r_\alpha}(x)=\mbox{Re}\pa{X e^{2i\pi x}}$,\\
	\item $\displaystyle\Lambda=\frac{1}{4}\E\int_{\T} \pa{\zetaup_1 '+(\overline{U}\bar{\zetaup}_1)'-(\overline{U}\bar{\zetaup}_1)'\circ r_\alpha}^2dx+O(\ep^3)=\frac{1}{8}\E\pa{|X|^2}+O(\ep^3).$
\end{itemize}
The result follows by simply rewriting $\E\pa{|X|^2}=\E\pa{|\overline{X}e^{2i\pi\alpha}|^2}$

\subsection{Proof of Theorem \ref{mainproj2}}
We are going to prove Theorem \ref{mainproj2} by mimicking the proof of Theorem \ref{principal}. Let $\delta>0$ and let $M$ be a random matrix in $SL_2(\R)$ such that $\|Tr(M)\|_{L^2(\Omega)}\leq 2-\delta$. Let $\alpha$ in $\T$ be so that $d(M,\mathcal{R})=||M-R_\alpha||$, and let $f_M=r_\alpha+\zetaup$ be the associated random diffeomorphism of $\T$. We assume that $M$ is valued in the open set 
$$\mathcal{U}_0=\{N\in Sl_2(\R), d(N,\mathcal{R}) < \beta\}$$ where $\beta$ is a constant depending only on $\delta$ and $\|\cdot\|$ chosen so that for $M$ in $\mathcal{U}_0$ we have  $|f_M'-1|\leq \frac{1}{2}$ and $|Tr(M)-Tr(R_\alpha)|\leq\frac{\delta}{2}$. The second inequality implies that $\|Tr(R_\alpha)|\|_{L^2(\Omega)}\geq 2-\frac{\delta}{2}$ and so $\|d(2\alpha,\Z)\|_{L^2(\Omega)}\geq \delta'$ for some positive $\delta'$ ($\approx \sqrt{\delta})$ depending on $\delta$, so the technics used to prove Theorem \ref{mainproj} still work.\\

Let us construct the first conjugation.
\begin{lem1}\label{firstproj}
	There exists $P$ in $SL_2(\R)$ such that either $||d(PMP^{-1},\mathcal{R})||_{L^2(\Omega)}\leq 4A_0\Lambda^{\frac{1}{2}}$ or $||d(PMP^{-1},\mathcal{R})||_{L^2(\Omega)}\leq C ||d(M,\mathcal{R})||_{L^2(\Omega)}^{\frac{3}{2}}$, where $A_0$ is the constant of Lemma \ref{equiv}, and $C$ is a constant depending only on $\delta$ and the norms. Moreover $\|P-I_2\|\leq C ||d(M,\mathcal{R})||_{L^2(\Omega)}$.
\end{lem1}

\begin{proof}
	From the proof of Lemma \ref{LyapuProj}, we get that setting $\etaup=\overline{U}\bar{\zetaup}_1$, $g=Id-\etaup$,  $\tilde{f}=gf_Mg^{-1}=r_\alpha+\tilde{\zetaup}$ and $\ep=||d(M,\mathcal{R})||_{L^2(\Omega)}$, we have
	$$\Lambda\geq \frac{1}{8}\int_{\T}\tilde{\zetaup}'^2 dx+O(\ep^3),$$
	using that if $\ep$ is small enough, $\tilde{f}'<2$ so $\ln(\tilde{f}')\leq\tilde{\zetaup}'-\frac{1}{4}\tilde{\zetaup}'^2$. So there exists 
	a constant $C$ such that
	$$\E\int_\T\tilde{\zetaup}'^2dx\leq 8\Lambda+C\ep^3,$$
	so 
	$$||d_0(\tilde{f},r_{\tilde{\alpha}})||_{L^2(\Omega)}\leq 3\Lambda^{\frac{1}{2}}+C^{\frac{1}{2}}\ep^{\frac{3}{2}}$$
	where $\tilde{\alpha}=\alpha+\int_{\T}\tilde{\zetaup}dx$.\\
	
	By Lemma \ref{projR}, there exists $P$ in $SL_2(\R)$ such that $\|P-I_2\|=O(\ep)$ and  $f_P(x)=x-\etaup(x)+O(||\etaup||^2)=g(x)+O(\ep^2)$. Let us set $\widetilde{M}=PMP^{-1}$. Since $d_0(f_P,g)=O(\ep^2)$, we deduce from Proposition \ref{esti0} that $d_0(f_{\widetilde{M}},\tilde{f})=d_0(f_Pf_Mf_P^{-1},gf_Mg^{-1})=O(\ep^2).$ Hence
	
	$$||d_0(f_{\widetilde{M}},r_{\tilde{\alpha}})||_{L^2(\Omega)}\leq 3\Lambda^{\frac{1}{2}}+C\ep^{\frac{3}{2}}$$
	for some new constant $C$. So either  $||d_0(f_{\widetilde{M}},r_{\tilde{\alpha}})||_{L^2(\Omega)}\leq 4\Lambda^{\frac{1}{2}}$ or $||d_0(f_{\widetilde{M}},r_{\tilde{\alpha}})||_{L^2(\Omega)}\leq 4C\ep^{\frac{3}{2}}$, and the conclusion follows from the inequality $||\widetilde{M}-R_{\tilde{\alpha}}||\leq A_0 d_0(f_{\widetilde{M}},r_{\tilde{\alpha}})$
	
\end{proof}

We can now prove Theorem \ref{mainproj2}.

\begin{proof}(Theorem \ref{mainproj2})\\
	Let $M$ be a random matrix with Lyapunov exponent $\Lambda$. We are going to assume that $d(M,\mathcal{R}) < \frac{\beta}{2}$ a.s. (in particular, $M\in \mathcal{U}_0$).  We construct a sequence of random matrices $(M_n)_n$ by induction: we set $M_0=M$, and once $M_n$ defined, if $||d(M_n,\mathcal{R})||_{L^2(\Omega)}\leq 4A_0\Lambda^{\frac{1}{2}}$ or if $M_n$ does not belong almost surely to $\mathcal{U}_0$ then we stop the sequence, and if not then we use Lemma \ref{firstproj} and set $M_{n+1}=P_n M_n P_n^{-1}$ where $P_n$ is given by the  lemma. Thus we get a sequence $(M_n)_{n\leq N}$ where $N$ belongs to $\N\cup\{+\infty\}$. Finally, we set $Q_n=P_{n-1}\cdots P_0$, so that $M_n=Q_{n}MQ_{n}^{-1}$.\\
	
	Let $\ep_n=||d(M_n,\mathcal{R})||_{L^2(\Omega)}$. By invariance by conjugation, the Lyapunov exponent of $M_n$ is $\Lambda$. So from the construction and Lemma \ref{firstproj} we deduce that for every $n<N$, $\ep_{n+1}\leq C\ep_n^{\frac{3}{2}}$ and for every $n\leq N$, $||P_n-I_2||\leq C\ep_n$. It is then straightforward that there is a constant $C_1$ and a positive number $\bar{\ep}$ such that if $\ep_0\leq \bar{\ep}$ then for every $n\leq N$, $\ep_n\leq C_1 2^{-\left(\frac{3}{2}\right)^n}\ep_0$,  and also $||Q_n-I_2||\leq C_1\ep_0$, and then that $d(M_n,\mathcal{R})\leq \beta $, i.e. $M_n\in\mathcal{U}_0$ (so the sequence only stop if $||d(M_n,\mathcal{R})||_{L^2(\Omega)}\leq 4A_0\Lambda^{\frac{1}{2}}$).\\
	
	Two cases can occur:
	\begin{itemize}
		\item If $\Gamma>0$, then $N<+\infty$. So $||d(M_N,\mathcal{R})||_{L^2(\Omega)}\leq 4A_0\Lambda^{\frac{1}{2}}$ with $M_N=Q_NMQ_N^{-1}$, and $||Q_N-I_2||\leq C_1\ep_0.$
		\item If $\Lambda=0$ then $N=+\infty$. Since $||Q_{n+1}-Q_{n}||=O(||Q_{n}||\cdot ||P_n-I_2||)=O(\ep_n)$, $(Q_n)$ converge to some matrix $Q$ such that $||Q-I_2||=O(\ep_0)$, and since $||d(Q_{n}MQ_{n}^{-1},\mathcal{R})||_{L^2(\Omega)}=\ep_n\to 0$, we conclude that $QMQ^{-1}\in \mathcal{R}$ almost surely.
	\end{itemize}
	Thorem \ref{mainproj2} follows.
\end{proof}

\section{Appendix: $C^k$ estimates} \label{appen}
In this section we give a quick proof of the propositions stated in Section \ref{prelim} and state some other classical $C^k$ estimates.\\

In the following propositions we consider maps $f:\R\rightarrow\R$. We denote by $\|\cdot\|_\infty$ the supremum norm, that is $\|f\|_\infty=\sup_\R |f|$.

\begin{prop1}(Kolmogorov inequality)\\
	For any integers $j\leq k$ and for any $f$ in $C^k(\R)$,
	$$\|f^{(j)}\|_\infty\leq C\|f^{(k)}\|_\infty^{j/k}\|f\|_\infty^{1-j/k}.$$
	where $C$ is a constant depending only on $k$.
\end{prop1}
\begin{proof}
	Being given real numbers $x$ and $h$, Taylor-Lagrange formula gives the existence of $c$ in $\R$ such that
	\begin{equation}\label{Taylor}f(x+h)=\sum_{n=0}^{k-1}f^{(n)}(x)\frac{h^n}{n!}+f^{(k)}(c)\frac{h^k}{k!}\end{equation}
	We fix real number $a_0,\ldots,a_{k-1}$ such that $\sum_{m=0}^{k-1}a_m n^m=\delta_{n,j}$ for $n=1,\ldots,k-1$ by inverting a Vandermonde system. For $t\in \R$ given, by a linear combinations of the formulas (\ref{Taylor}) with $h=0,t,2t,\ldots, (k-1)t$ we get
	$$\sum_{m=0}^{k-1}a_m f(x+mt)=f^{(j)}(x)\frac{t^j}{j!}+\left(\sum_{m=0}^{n-1}a_mf^{(k)}(c_m)\right)\frac{t^k}{k!}$$
	for some real numbers $c_1,\ldots,c_{k-1}$. In particular,
	$$\|f^{(j)}\|_\infty\leq C(t^{-j}\|f\|_\infty+t^{k-j}\|f^{(k)}\|_\infty)$$
	for some constant $C$, and the result follows by optimizing in $t$.	
\end{proof}
\begin{prop1}\label{Kolmocor}(Product of norms of derivatives)\\
	For any $f$, $g$ in $C^k(\R)$, and any integer $j\leq k$,
	$$\|f^{(j)}\|_\infty\|g^{(k-j)}\|_\infty\leq C(\|f^{(k)}\|_\infty\|g\|_\infty+\|f\|_\infty\|g^{(k)}\|_\infty)$$
	where $C$ is a constant depending only on $k$.
\end{prop1}
\begin{proof}
	It is a consequence of Kolmogorov inequality and the convexity inequality $a^{\theta}b^{1-\theta}\leq \theta a+(1-\theta)b$:
	$$\|f^{(j)}\|_\infty\|g^{(k-j)}\|_\infty\leq C\|f^{(k)}\|_\infty^{j/k}\|f\|_\infty^{1-j/k}\|g^{(k)}\|_\infty^{1-j/k}\|g\|_\infty^{j/k}\leq C\left(\frac{j}{k}\|f^{(k)}\|_\infty\|g\|_\infty+(1-\frac{j}{k})\|f\|_\infty\|g^{(k)}\|_\infty\right).$$
\end{proof}

\begin{prop1}(Derivative of a product)\label{estiprod}
	For any integer $k$ and any $f$, $g$ in $C^k(\R)$,
	$$\|(fg)^{(k)}\|_\infty\leq C(\|f^{(k)}\|_\infty\|g\|_\infty+\|f\|_\infty\|g^{(k)}\|_\infty)$$
	where $C$ is a constant depending only on $k$.
\end{prop1}
\begin{proof}
	By Leibnitz formula, $\|(fg)^{(k)}\|_\infty\leq \sum_{j=0}^k \begin{pmatrix}
	k\\j
	\end{pmatrix} \|f^{(j)}\|_\infty \|g^{(k-j)}\|_\infty$, and by the proposition above, $\|f^{(j)}\|_\infty\|g^{(k-j)}\|_\infty\leq C(\|f^{(k)}\|_\infty\|g\|_\infty+\|f\|_\infty\|g^{(k)}\|_\infty)$ for some $C$.
\end{proof}
\begin{prop1}(Derivative of a composition)\label{esticomp}
	Let $M\geq 1$. For any integer $k\geq 1$ and any $f$, $g$ in $C^k(\R)$ such that $|g'|\leq M$ on $\R$,
	$$\|(f\circ g)^{(k)}\|_\infty\leq CM^{k-1}(\|f^{(k)}\|_\infty\|g'\|_\infty+\|f'\|_\infty\|g^{(k)}\|_\infty)$$
	where $C$ is a constant depending only on $k$.
\end{prop1}
\begin{proof}
	We proceed by induction on $k$. The statement is obvious for $k=1$. Let $k\geq 2$. Since $(f\circ g)^{(k)}=\left(f'\circ g\cdot g'\right)^{(k-1)}$, we obtain by Proposition \ref{estiprod} for some constant $C$
	$$\|(f\circ g)^{(k)}\|_\infty \leq C\left(\|(f'\circ g)^{(k-1)}\|_\infty\|g'\|_\infty+\|f'\circ g\|_\infty\|(g')^{(k-1)}\|_\infty\right),$$
	so
	$$\|(f\circ g)^{(k)}\|_\infty\leq C\left(M\|(f'\circ g)^{(k-1)}\|_\infty+\|f'\|_\infty\|g^{(k)}\|_\infty\right).$$
	By induction hypothesis,
	$$\|(f'\circ g)^{(k-1)}\|_\infty\leq CM^{k-2}(\|f^{(k)}\|_\infty\|g'\|_\infty+\|f''\|_\infty\|g^{(k-1)}\|_\infty)$$
	for some constant $C$ depending on $k$. So for some new constant $C$
	$$\|(f\circ g)^{(k)}\|_\infty\leq CM^{k-1}\left(\|f^{(k)}\|_\infty\|g'\|_\infty+\|f''\|_\infty\|g^{(k-1)}\|_\infty+\|f'\|_\infty\|g^{(k)}\|_\infty\right).$$
	By Proposition \ref{Kolmocor}
	$$\|f''\|_\infty\|g^{(k-1)}\|_\infty \leq C(\|f^{(k)}\|_\infty\|g'\|_\infty+\|f'\|_\infty\|g^{(k)}\|_\infty)$$
	for some constant $C$ so finally, with a new constant $C$,
	$$\|(f\circ g)^{(k)}\|_\infty\leq CM^{k-1}(\|f^{(k)}\|_\infty\|g'\|_\infty+\|f'\|_\infty\|g^{(k)}\|_\infty),$$
	which completes the induction.
\end{proof}

From these general estimates, we deduce some more specific ones for our context. We reintroduce the $C^k$-norms: for $\phi$ in $C^k(\R)$ we define its $C^k$-norm by $\|\phi\|_k=\max(\|\phi\|_{\infty},\|\phi'\|_{\infty},\ldots,\|\phi^{(k)}\|_{\infty})$ (in particular, $\|\cdot\|_0$ is also the supremum norm). Alternatively we could define $\|\phi\|_k=\max(\|\phi\|_{\infty},\|\phi^{(k)}\|_{\infty})$, which is an equivalent norm thanks to Kolmogorov inequality.

\begin{lem1}\label{compodif}
	Let $k$ be an integer, let $M\geq 1$, let $f$, $g$ be in $C^k(\R)$ such that $|f'|, |g'|\leq M$ on $\R$. Then:
	
	$$\|f\circ g-Id\|_k\leq CM^k(\|f-Id\|_k+\|g-Id\|_k)$$
	where $C$ is a constant depending only on $k$.
\end{lem1}

\begin{proof}
	Let $\varphi=f-Id$ and $\psi=g-Id$. Since $f\circ g-Id=\psi+\varphi\circ g$,
	we only need to bound $\|\varphi\circ g\|_k$. We have $\|\varphi\circ g\|_0=\|\varphi\|_0$, $\|(\varphi\circ g)'\|_0\leq \|g'\|_0\|\varphi'\|_0\leq M\|\varphi\|_1$, and if $k\geq 2$, by Proposition \ref{esticomp} for some constant $C$ depending on $k$ we have
	$$\|(\varphi\circ g)^{(k)}\|_0\leq CM^{k-1}(\|\varphi^{(k)}\|_0\|g'\|_0+\|\varphi'\|_0\|g^{(k)}\|_0),$$
	with $\|\varphi'\|_0 \leq 1+M\leq 2M$, $\|g'\|_0\leq M$ and $\|g^{(k)}\|_0=\|\psi^{(k)}\|_0$, so 
	\begin{equation}\label{comp}\|\varphi\circ g\|_k\leq CM^k(\|\varphi\|_k+\|\psi\|_k),\end{equation}
	for some new constant $C$ depending on $k$, and the statement follows.
\end{proof}

\begin{lem1}\label{invdif}
	Let $k$ be an integer, let $q<1$, and let $f$ be in $C^k(\R)$ such that $|f'-1|\leq \frac{1}{2}$ on $\R$. Then :
	$$\|f^{-1}-Id\|_k\leq C\|f-Id\|_k$$
	where $C$ is a constant depending only on $k$.
\end{lem1}

\begin{proof}
	Let $g=f^{-1}$, $\varphi=f-Id$ and $\psi=g-Id$, so that the identity $f\circ g=Id$ becomes $\psi=-\varphi\circ g$. We want to prove that $\|\psi\|_k\leq C\|\varphi\|_k$ for some constant $C$. It is straightforward if $k=0$ or $1$ so we assume tha $k\geq 2$ and we make the induction assumption that for every $j<k$, $\|\psi\|_j\leq C\|\varphi\|_j$ for some constant $C$. Then:
	$$\|\psi\|_k=\|\varphi\circ g\|_k\leq\|\varphi\|_0+\|\varphi '\circ g\cdot g '\|_{k-1}\leq\|\varphi\|_0+\sum_{j=0}^{k-1}  \left( \begin{array}{c}k-1 \\ j\end{array}\right) \|\varphi '\circ g\|_{j}\|g '\|_{k-1-j}.$$
	For $j=0$, 
	$$\|\varphi '\circ g\|_{0}\|g '\|_{k-1}\leq \|\varphi'\|_{0}(1+\|\psi '\|_{k-1}) \leq \|\varphi\|_1+\frac{1}{2}\|\psi\|_k,$$ 
	and for $j\not=0$, by using inequality (\ref{comp}) (with $M=2$) and the induction assumption we can bound $\|\varphi '\circ g\|_{j}\leq C\|\varphi\|_j$ for some constant $C$, and then by using Proposition \ref{Kolmocor} we get  $\|\varphi '\circ g\|_{j}\|g '\|_{k-1-j}\leq C\|\varphi\|_k$ with a new constant $C$. So we deduce finally that we have for some constant $C$
	$$\|\psi\|_k\leq \frac{1}{2}\|\psi\|_k+C\|\varphi\|_k,$$ 
	and so $\|\psi\|_k\leq 2C\|\varphi\|_k$, which completes the induction.
\end{proof}

\begin{lem1}(a $C^k$ mean value inequality) \label{meanvalue}
	Let $M\geq 1$, let $f$, $g$  be in $C^k(\R)$ such that $|f'|,|g'|,|f^{(k)},|g^{(k)}|\leq M$ on $\R$, and let $\phi\in C^{k+1}(\R)$.
	Then 
	$$\|\phi\circ f-\phi\circ g\|_k\leq C\|\phi\|_{k+1}\|f-g\|_k$$
	where $C$ depends only on $k$ and $M$.
\end{lem1}
\begin{proof}
	We write 
	$$\phi\circ f-\phi\circ g=(f-g)\int_0^1\phi'\circ h_t dt$$
	where $h_t=(1-t)f+tg$. Thus,
	$$\|\phi\circ f-\phi\circ g\|_k\leq C\|f-g\|_k\int_0^1\|\phi'\circ h_t\|_k dt$$
	for some constant $C$ depending only on $k$. By Proposition \ref{esticomp} (and Kolmogorov inequality), $\|\phi'\circ h_t\|_k\leq  C\|\phi\|_{k+1}$ for some constant $C$ depending on $k$ and $M$. The result follows.
\end{proof}

Finally, let us prove Propositions \ref{estiK}, \ref{esti0}, \ref{esti1}, \ref{estival} of Section \ref{prelim}. Proposition \ref{estiK} is an immediate consequence of Lemmas \ref{compodif} and \ref{invdif} and the fact that $d_k$ is invariant by (left or right) composition by rotations. Propostion \ref{estival} is a straightforward consequence of inequality (\ref{comp}) since $d_k(f\circ h, g\circ h)=\|(f-g)\circ h\|_k$. To prove Propostion \ref{esti1}, we write $f=r_\alpha+\zetaup$ and $g=Id+\etaup$, and then an algebraic computation
gives
$$g\circ f \circ g^{-1}=r_\alpha+\left(\zetaup\circ g^{-1}+\eta\circ (f\circ g^{-1})-\eta\circ g^{-1}\right).$$
The difference between this map and the approximation $r_\alpha+\left(\zetaup+\eta\circ r_\alpha-\eta\right)$ can be estimated in $C^1$-norm thanks to Lemma \ref{meanvalue} (with $k=1$), what gives the result (alternatively one can directly bound this difference and its derivative by elementary calculus). Finally, Proposition \ref{esti0} is an elementary consequence of the invariance of $d_0$ by right composition and mean value inequality:
$$\begin{array}{ll}d_0(gfg^{-1},\tilde{g}f\tilde{g}^{-1})&\leq d_0(gfg^{-1},\tilde{g}fg^{-1})+d_0(\tilde{g}f\tilde{g}^{-1},\tilde{g}fg^{-1})\\
&\leq d_0(g,\tilde{g})+d_0(\tilde{g}f\tilde{g}^{-1}g,\tilde{g}f)\\
&\leq d_0(g,\tilde{g})+d_0((\tilde{g}f\tilde{g}^{-1})\circ g,(\tilde{g}f\tilde{g}^{-1})\circ \tilde{g})\\
&\leq (1+\|(\tilde{g}f\tilde{g}^{-1})'\|_0)d_0(g,\tilde{g}),
\end{array}$$
with $\|(\tilde{g}f\tilde{g}^{-1})'\|_0$ easily bounded by above.
\bibliographystyle{plain}
\bibliography{Bibliography}
\end{document}